\documentclass[12pt,a4paper,english,british,times,doublespace]{qjrms4}
\usepackage{mathptmx}
\usepackage[T1]{fontenc}
\usepackage[latin9]{inputenc}
\usepackage{url}
\usepackage{amsmath}
\usepackage{graphicx}
\usepackage{esint}
\usepackage[authoryear]{natbib}

\makeatletter

\pdfpageheight\paperheight
\pdfpagewidth\paperwidth

\providecommand{\tabularnewline}{\\}

\usepackage { fancybox}
\usepackage[export]{adjustbox}

\usepackage{todonotes}

\usepackage{amssymb}

\usepackage{lineno}

\makeatother

\usepackage{babel}
\begin{document}
\title{Dimension Splitting and a Long Time-Step Multi-Dimensional Scheme for Atmospheric Transport}
\author{Yumeng Chen,\affil{a} Hilary Weller,\affil{b}\corrauth\ Stephen Pring\affil{c} and James Shaw\affil{b}}

\runningheads{Y. Chen, H. Weller, S. Pring and J. Shaw}{Dimension splitting for Advection}

\address{
\affilnum{a}Universit\"at Hamburg\\
\affilnum{b}Meteorology, University of Reading \\
\affilnum{c}Met Office, UK}
\corraddr{E-mail: <h.weller@reading.ac.uk>}

\begin{abstract}\normalfont
Dimensionally split advection schemes are attractive for atmospheric modelling due to their efficiency and accuracy in each spatial dimension. Accurate long time steps can be achieved without significant cost using the flux-form semi-Lagrangian technique. The dimensionally split scheme used in this paper is constructed from the one-dimensional Piecewise Parabolic Method and extended to two dimensions using COSMIC splitting. The dimensionally split scheme is compared with a genuinely multi-dimensional, method of lines scheme with implicit time-stepping which is stable for very large Courant numbers.

Two-dimensional advection test cases on Cartesian planes are proposed that avoid the complexities of a spherical domain or multi-panel meshes. These are solid body rotation, horizontal advection over orography and deformational flow. The test cases use distorted meshes either to represent sloping terrain or to mimic the distortions of a cubed sphere.

Such mesh distortions are expected to accentuate the errors associated with dimension splitting, however, the dimensionally split scheme is very accurate on orthogonal meshes and accuracy decreases only a little in the presence of large mesh distortions. The dimensionally split scheme also loses some accuracy when long time-steps are used. The multi-dimensional scheme is almost entirely insensitive to mesh distortions and asymptotes to second-order accuracy at high resolution. As is expected for implicit time-stepping, phase errors occur when using long time-steps but the spatially well resolved features are advected at the correct speed and the multi-dimensional scheme is always stable.

A naive estimate of computational cost (number of multiplies) reveals that the implicit scheme is  the most expensive, particularly for large Courant numbers. If the multi-dimensional scheme is used instead with explicit time-stepping, the Courant number is restricted to less than one but the cost becomes similar to the dimensionally split scheme.
\end{abstract}

\keywords{
Multi-dimensional; advection; stability; accuracy; long time-steps
}
\maketitle

\onecolumn

\section{Introduction}

Weather and climate models are being developed on quasi-uniform meshes
in order to better exploit modern computers \citep[eg ][]{WTC12,LUJ+14,ST12,SG11,KNK15}
and so accurate and efficient transport (or advection) schemes on
non-orthogonal meshes are required. There is an abundance of desirable
properties of advection schemes, including:
\begin{enumerate}
\item Inherent local conservation of the advected quantity.
\item Stability in the presence of large Courant numbers.
\item Accuracy in the presence of large Courant numbers.
\item High order-accuracy.
\item Low computational cost, good parallel scaling and multi-tracer efficiency.
\item Low phase and dispersion errors (advection of all wavenumbers of the
advected quantity at close to the correct speed).
\item Low diffusion errors (maintaining amplitude of all wavenumbers of
the advected quantity).
\item Boundedness, monotonicity, positivity and maintaining correlations
between multiple advected tracers.
\end{enumerate}
Four (important) properties are listed together in the final item
because they will not be addressed here. In this paper, we address
the issue of conservative, accurate, efficient advection schemes for
logically rectangular, non-orthogonal meshes which are stable in the
presence of large Courant numbers. These schemes would be particularly
relevant for cubed-sphere meshes and for terrain following meshes.
Another novel aspect of this paper is that, for simplicity, we test
advection schemes entirely on planar meshes rather than on the sphere,
proposing test cases to challenge advection schemes on non-orthogonal
meshes without the need to implement meshes in spherical geometry.

Dimensionally split schemes (operating separately in each spatial
dimension) are attractive for atmospheric modelling due to their efficiency
and high accuracy in each spatial dimension \citep[eg. ][]{LR96,LLM96,BS99,PL07,KNK15}.
Inherent conservation is guaranteed by using the flux-form semi-Lagrangian
(FFSL or forward in time) technique \citep[eg. ][]{CW84} which integrates
the dependent variable over a swept distance upstream of every face
in order to calculate fluxes in and out of cells. Accurate long time
steps can be achieved without significant cost by calculating cumulative
mass fluxes along the domain \foreignlanguage{english}{\citep{LLM95}.
}This is done without remapping which can be an expensive procedure,
finding a conservative map between fields on overlapping meshes\foreignlanguage{english}{.
}One-dimensional schemes can be used with operator splitting to create
dimensionally split, second-order accurate schemes \citep[eg][]{LLM96}
on logically rectangular, multidimensional meshes. Dimensionally split
schemes have been found to give good accuracy on non-orthogonal meshes
such as the cubed-sphere \citep{PL07,KNK15} with special treatment
over cube edges. \citet{PL07} use the average of two one-sided schemes
at cube edges whereas \citet{KNK15} create ghost cells outside each
cube panel boundary. Without this special treatment, dimensionally
split schemes do not account properly for mesh distortion and errors
may occur. Accuracy between second and fourth order was found in practice
on a range of test cases on the cubed-sphere by \citet{KNK15}. 

The same problems occur when using dimensionally split schemes over
orography since terrain following layers become non-orthogonal over
orography and special treatment cannot be used everywhere where there
is orography. A common solution has been to make the terrain following
layers as smooth as possible, reducing non-orthogonality \citep[eg][]{SLF+02}
or to use floating Lagrangian vertical co-ordinates \citep{Lin04}.
Dimension splitting may account for some of the errors over orography
reported by \citet{KUJ14} although they do not cite this as a reason
for errors.

Dimension splitting errors on distorted meshes can be eliminated by
using genuinely multi-dimensional advection schemes. These can be
either FFSL \citep[ie swept area, eg.][]{Las02,LR05b,Miu07b,TCD14},
method of lines \citep[MOL, discretising space and time separetely, eg.][]{WWF09,SG11,KNK15,SWMD1x}
or conservative semi-Lagrangian \citep[with conservative re-mapping, eg][]{IK04,LNU10,ZWS04}.
The FFSL and MOL multi-dimensional schemes have not previously been
extended to work with Courant numbers significantly larger than one.
FFSL multi-dimensional schemes could be extended to handle large time-steps
by integrating the upstream swept volume over a large upstream volume,
interacting with a large number of upstream cells. However the cost
would be proportional to the time-step since the larger the Courant
number, the more cells the upstream swept volume would need to overlap
with. This technique therefore offers no advantage over using a smaller
time-step. MOL multi-dimensional schemes can be extended to work with
Courant numbers larger than one by using implicit time-stepping. This
will increase the computational cost per tracer advected since the
solution of a matrix equation would be needed for every advected tracer.
Two other disadvantages of implicit time-stepping are the large phase
errors when long time-steps are used \citep[eg][]{DB12,LWW14} and
the difficulty of achieving monotonicity. This is in contrast to semi-Lagrangian
or FFSL schemes which maintain accuracy with long time-steps \citep{Pur76,PS84,LLM95}
although monotonocity with long time-steps is still challenging \citep{Bott10}.
Conservative semi-Lagrangian naturally extends to long time-steps
but the conservative remapping is complicated and expensive, particularly
on non-rectangular meshes and will not be investigated here. \citet{LUJ+14}
described how the FFSL technique with a long time-step can be made
equivalent to the conservative semi-Lagrangian.

It is therefore not clear what approach should be taken for achieving
long time-steps when advecting multiple tracers on distorted meshes.
In this paper, we show the effect of dimension splitting errors using
a FFSL dimensionally split scheme on a number of test-cases which
use distorted meshes and compare with a genuinely multi-dimensional
implicit MOL scheme using large and small Courant numbers.

The theoretical properties of dimensionally split advection schemes
are often tested on uniform, orthogonal meshes \citep[eg][]{LLM96}.
On the cubed-sphere, special treatment is needed for the cube edges
\citep[eg][]{LR96,KNK15}. Developing a transport scheme to the extent
that it can be used on a multi-panel cubed-sphere with special treatment
of cube edges is a considerable undertaking. Hence there is a need
for more challenging advection test cases which are simpler to implement,
without the need for spherical meshes. We therefore propose some modifications
of existing test cases to use distorted meshes, or distorted co-ordinate
systems, on a logically rectangular, two-dimensional plane. 

The long time-step permitting, dimensionally split scheme and the
long time-step permitting multi-dimensional scheme are defined in
section \ref{sec:models}. In section \ref{sec:results} we present
results of three advection test cases on distorted meshes in two-dimensional
planes using Courant numbers above and below one. These are the solid
body rotation test case of \citet{LLM96} modified to use a mesh (or
co-ordinate system) with distortions similar to a cubed-sphere (section
\ref{sub:solidRotate}), the horizontal advection test case over orography
\citep{SLF+02}, examining sensitivity to time-step, resolution and
mountain height, all on the maximally distorted basic terrain following
mesh (section \ref{sub:overOrog}) and a modification of the deformational
flow test case of \citet{LSPT12} for a periodic rectangular plane
(section \ref{sub:deformFlow}). Some estimates are made of computational
cost in section \ref{sub:cost} and final conclusions are drawn in
section \ref{sec:concs}.

\section{Transport Schemes\label{sec:models}}

We present two conservative advection schemes suitable for long time-steps
(stable for Courant numbers significantly larger than one) for solving
the linear advection equation:
\begin{equation}
\frac{\partial\phi}{\partial t}+\nabla\cdot\mathbf{u}\phi=0
\end{equation}
where the dependent variable $\phi$ is advected by velocity field
$\mathbf{u}(\mathbf{x},t)$. The dimensionally split scheme is the\foreignlanguage{english}{
piecewise parabolic method \citep[PPM][]{CW84}} which uses the flux-form
semi-Lagrangian approach extended to long time-steps following \citet{LLM95}
with COSMIC splitting \citep{LLM96} to extend PPM to two dimensions.
The multi-dimensional scheme uses the method of lines approach (treating
space and time independently). The second-order accurate spatial discretistaion
of \citet{WS14} is combined with implicit, Crank-Nicholson time-stepping
to allow long time-steps. Neither scheme has monotonicity or positivity
preservation. 

The code for dimensionally splitting scheme is available at \url{https://github.com/yumengch/COSMIC-splitting}.
The multi-dimensional scheme is implemented using \citet{OpenFOAM16}
and is available at \url{https://github.com/AtmosFOAM/}.

\subsection{One-dimensional PPM with Long Time-steps\label{sub:PPM}}

We describe a long time-step version of PPM\foreignlanguage{english}{
for solving the one-dimensional advection equation:
\begin{equation}
\frac{\partial\phi}{\partial t}+\frac{\partial u\phi}{\partial x}=0.
\end{equation}
\citet{CW84} defined PPM with monotonicity constraints and for variable
resolution but for simplicity (and for comparison with the multi-dimensional
scheme) we will define PPM without monotonicity constraints and for
a fixed resolution, $\Delta x$, and time-step, $\Delta t$. With
these restrictions, PPM should be fourth-order accurate in one dimension
for vanishing time-step. We define the dependent variable, $\phi_{i}^{(n)}$,
to be the mean value of $\phi$ in cell $i$ at time-level $n$ where
$x_{i}=i\Delta x$ and $t=n\Delta t$. Since PPM is a flux-form finite
volume method, $\phi_{i}^{(n)}$ is updated using:
\begin{equation}
\phi_{i}^{(n+1)}=\phi_{i}^{(n)}+X_{C}\left(\phi\right)=\phi_{i}^{(n)}-\frac{u_{i+1/2}\phi_{i+1/2}-u_{i-1/2}\phi_{i-1/2}}{\Delta x}\label{PPMgauss}
\end{equation}
where $X_{C}\left(\phi\right)$ is the conservative advection operator
for $\phi$ in the $x$ direction. The fluxes, $u_{i\pm1/2}\phi_{i\pm1/2}$,
are found by integrating a piecewise polynomial, $p$, along the distance
travelled in each time-step upwind of cell boundary $x_{i\pm1/2}$.
The polynomial is defined in each cell, $i$, such that:
\begin{equation}
\phi_{i}=\frac{1}{\Delta x}\int_{x_{i-1/2}}^{x_{i+1/2}}p_{i}(x)\ dx
\end{equation}
by
\begin{equation}
p_{i}(x)=p_{i-1/2}+\xi\left(p_{i+1/2}-p_{i-1/2}+\left(1-\xi\right)6\left(\phi_{i}-\frac{1}{2}\left(p_{i-1/2}+p_{i+1/2}\right)\right)\right)
\end{equation}
where $\xi=(x-x_{i-1/2})/\Delta x$ and 
\begin{equation}
p_{i+1/2}=\frac{7}{12}(\phi_{i}+\phi_{i+1})-\frac{1}{12}(\phi_{i+2}+\phi_{i-1})
\end{equation}
}

In order to cope with long time-steps, we follow \citet{LLM95} and
divide the Courant number into a signed integer part, $c_{N}$, and
a remainder, $c_{r}$. The departure point of location $x_{i-\frac{1}{2}}$
is thus computed as $x_{d}=x_{i-\frac{1}{2}}-u_{i-\frac{1}{2}}\Delta t$
and the departure cell, for cell edge $i-\frac{1}{2}$, is $i_{d}=i-c_{N}-1$
for $u_{i-\frac{1}{2}}>0$ and $i_{d}=i-c_{N}$ for $u_{i-\frac{1}{2}}<0$.
Then for $u_{i-\frac{1}{2}}>0$ the flux through $x_{i-1/2}$ between
times $n\Delta t$ and $(n+1)\Delta t$ is: 
\begin{equation}
u_{i-\frac{1}{2}}\phi_{i-\frac{1}{2}}=\frac{1}{\Delta t}\int_{x_{d}}^{x_{i-1/2}}p(x)\ dx=\frac{1}{\Delta t}\left(M_{i-1/2}-M_{i-c_{N}-1/2}+\int_{x_{i-1/2-c_{N}}-c_{r}\Delta x}^{x_{i-1/2-c_{N}}}p(x)\ dx\right)\label{eq:PPMflux}
\end{equation}
where $M_{i-1/2}$ is the cumulative mass from the start point to
position $x_{i-1/2}$: 
\begin{equation}
M_{i-1/2}=\sum_{k<i}\Delta x\phi_{k}.\label{eq:cumulativeMass}
\end{equation}
This departure point calculation assumes that the velocity is uniform
on the computational mesh which has a first-order error which could
be particularly damaging for long Courant numbers, when the wrong
departure cell could be found.

The velocity is derived from a stream function and the Jacobian of
the co-ordinate transform: \foreignlanguage{british}{
\begin{equation}
\left(\begin{array}{c}
u\\
v
\end{array}\right)=J\left(\begin{array}{c}
\;\Psi_{y}\\
-\Psi_{x}
\end{array}\right)
\end{equation}
}For stability, the time-step is restricted by the deformational Courant
number:
\begin{equation}
c_{d}=\Delta t\max\left(\biggl|\frac{\partial u}{\partial x}\biggr|,\biggl|\frac{\partial u}{\partial y}\biggr|,\biggl|\frac{\partial u}{\partial z}\biggr|,\biggl|\frac{\partial v}{\partial x}\biggr|,\biggl|\frac{\partial v}{\partial y}\biggr|,\biggl|\frac{\partial v}{\partial z}\biggr|,\biggl|\frac{\partial w}{\partial x}\biggr|,\biggl|\frac{\partial w}{\partial y}\biggr|,\biggl|\frac{\partial w}{\partial z}\biggr|\right)\label{eq:cd}
\end{equation}
\citep{PS84} such that $c_{d}\leq1$.

\subsection{COSMIC Splitting\label{sub:COSMIC}}

COSMIC operator splitting \citep{LLM96} allows single stage, one-dimensional
schemes such as PPM to be used stably in two or more dimensions whilst
retaining conservation, constancy preservation and second-order accuracy
(on orthogonal meshes). As we are now considering two spatial dimensions,
we define $\phi_{ij}$, $u_{ij}$ and $v_{ij}$, the values of $\phi$
and the velocity components, $u$ and $v$ in cell $(i,j)$ where
$x=i\Delta x$ and $y=j\Delta y$. COSMIC splitting uses both advective
and conservative advection operators in the $x$ and $y$ directions:
\begin{align}
X_{C}(\phi)=-\frac{1}{\Delta x}\left(u_{e}\phi_{e}-u_{w}\phi_{w}\right) &  & Y_{C}(\phi)=-\frac{1}{\Delta y}\left(v_{n}\phi_{n}-v_{s}\phi_{s}\right)\label{eq:COSMIC_A}\\
X_{A}(\phi)=X_{c}(\phi)+\frac{\phi_{ij}}{\Delta x}\left(u_{e}-u_{w}\right) &  & Y_{A}(\phi)=Y_{c}(\phi)+\frac{\phi_{ij}}{\Delta y}\left(v_{n}-v_{s}\right)\label{eq:COSMIC_C}
\end{align}
where $\phi_{n}=\phi_{i,j+1/2}$, $\phi_{s}=\phi_{i,j-1/2}$, $\phi_{e}=\phi_{i+1/2,j}$
, $\phi_{w}=\phi_{i-1/2,j}$, $v_{n}=v_{i,j+1/2}$, $v_{s}=v_{i,j-1/2}$,
$u_{e}=u_{i+1/2,j}$ and $u_{w}=u_{i-1/2,j}$ are the values of $\phi$,
$u$ and $v$ at the cell boundaries. If COSMIC is being used to extend
PPM to two spatial dimensions then $\phi_{n,s,e,w}$ are calculated
from equation \ref{eq:PPMflux}. Assuming C-grid staggering, $v_{n}$,
$v_{s}$, $u_{e}$ and $u_{w}$ are dependent variables. Instead of
using cell centered velocity \citep{LR96} or upwind velocity \citep{LLM96},
the advective operators are calculated in a similar manner to \citet{Lin04}.

Mesh distortions can be included in the advection equation with a
co-ordinate transform with Jacobian $J$:\foreignlanguage{english}{
\begin{equation}
\frac{\partial|J|^{-1}\phi}{\partial t}+\frac{\partial|J|^{-1}u\phi}{\partial x}+\frac{\partial|J|^{-1}v\phi}{\partial y}=0.
\end{equation}
}The operators $X_{c}$, $Y_{c}$ $X_{A}$ and $Y_{A}$ are then:\foreignlanguage{english}{
\begin{align}
X_{C}(\phi)=-\frac{1}{\Delta x}\left(|J|_{e}^{-1}u_{e}\phi_{e}-|J|_{w}^{-1}u_{w}\phi_{w}\right) &  & Y_{C}(\phi)=-\frac{1}{\Delta y}\left(|J|_{n}^{-1}v_{n}\phi_{n}-|J|_{s}^{-1}v_{s}\phi_{s}\right)\label{eq:COSMIC_A-1}\\
X_{A}(\phi)=X_{c}(\phi)+\frac{\phi_{ij}}{\Delta x}\left(|J|_{e}^{-1}u_{e}-|J|_{w}^{-1}u_{w}\right) &  & Y_{A}(\phi)=Y_{c}(\phi)+\frac{\phi_{ij}}{\Delta y}\left(|J|_{n}^{-1}v_{n}-|J|_{s}^{-1}v_{s}\right)\label{eq:COSMIC_C-1}
\end{align}
which are combined to update $\phi_{ij}$ in each cell by: 
\begin{equation}
\phi_{ij}^{(n+1)}=\phi_{ij}^{n}+|J|_{ij}X_{C}\left(\phi_{ij}^{(n)}+\frac{|J|_{ij}}{2}Y_{A}\left(\phi_{ij}^{(n)}\right)\right)+|J|_{ij}Y_{C}\left(\phi_{ij}^{(n)}+\frac{|J|_{ij}}{2}X_{A}\left(\phi_{ij}^{(n)}\right)\right).\label{eq:COSMICupdate}
\end{equation}
where $|J|_{ij}$ is the determinant of Jacobian at the center of
cell $(i,j)$, and $\Delta x$ and $\Delta y$ give the cell size
in uniform computational domain.}

\subsection{Multi-dimensional Method of Lines (MOL) Scheme\label{sub:MMmolscheme}}

\begin{figure}
\noindent \begin{centering}
\includegraphics{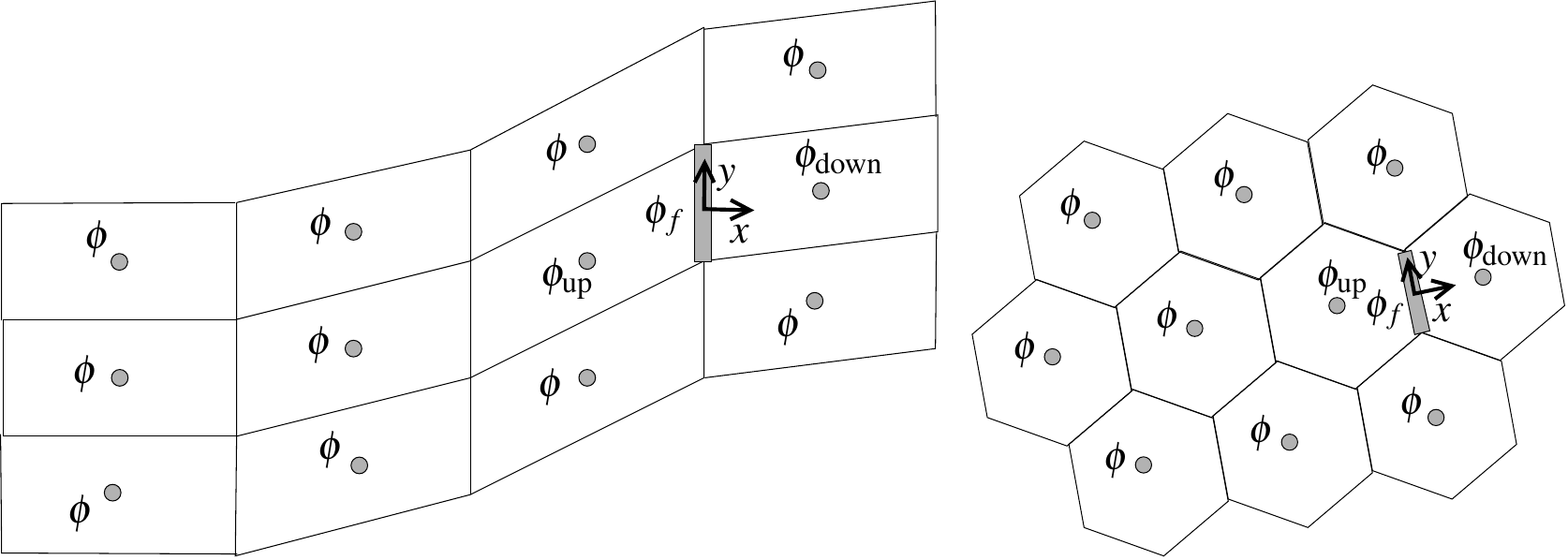}
\par\end{centering}

\caption{Stencils of upwind-biased cells for interpolating onto face $f$ using
a 2d cubic polynomial for two different mesh structures.\label{fig:stencil}}
\end{figure}

The MOL scheme uses the finite volume method for arbitrary meshes
and is implemented in \citep{OpenFOAM16}. This uses a cubic upwind
spatial discretisation \citep{WS14,SWMD1x} combined with Crank-Nicholson
in time. Although the interpolation uses a cubic polynomial, cell
centre values are approximated as cell average values and face centre
values are approximated as face average values so the method is limited
to 2nd order accuracy in space. The advection scheme uses Gauss's
divergence theorem to approximate the divergence term of the advection
equation:
\begin{equation}
\nabla\cdot\mathbf{u}\phi\approx\frac{1}{V}\sum_{f\in c}\phi_{f}\mathbf{u}_{f}\cdot\mathbf{S}_{f}\label{eq:multidGauss}
\end{equation}
where $V$ is the cell volume, the summation is over all faces, $f$,
of cell $c$, $\phi_{f}$ is the value of the dependent variable,
$\phi$ interpolated onto face $f$, $\mathbf{u}_{f}$ is the velocity
at face $f$ and $\mathbf{S}_{f}$ is the vector normal to face $f$
with magnitude equal to the area of face $f$ (ie the face area vector).
The advection scheme uses a fit to a 2d (or 3d) polynomial using an
upwind-biased stencil of cells (fig \ref{fig:stencil}) in order to
interpolate from known cell values onto a face. In 2d, the cubic polynomial
is:\foreignlanguage{english}{
\begin{equation}
\phi=a+bx+cy+dx^{2}+exy+fy^{2}+gx^{3}+hx^{2}y+ixy^{2}
\end{equation}
omitting the $y^{3}$ term, where $x$ is the direction normal to
a cell face and $y$ is perpendicular to $x$. (The $y^{3}$ term
is omitted because it cannot be set with a stencil that is narrow
in the direction of the flow.) Coefficients $a$ to $i$ are set from
a least squares fit to the cell data in the stencil. The least-squares
problem involves a $9\times m$ matrix singular value decomposition
(where $m$ is the size of the stencil) for every face and for both
orientations of each face . However this is purely a geometric calculation
and is therefore a pre-processing activity since the mesh is fixed.
This generates a set of weights for calculating $\phi_{f}$ from the
cell values in the stencil, leaving $m$ multiplies for each face
for each call of the advection operator. The stencils are found for
three-dimensional, arbitrarily structured meshes by finding the face(s)
closest to upwind of the face we are interpolating onto, taking the
two cells either side of the upwind face(s) and then taking the vertex
neighbours of those central cells (fig \ref{fig:stencil}). For each
face there are two possible stencils depending on the upwind direction.
Both of the stencils are stored and the interpolation weights for
both stencils are calculated. }

In order to ensure that the fit is accurate in the cells either side
of face $f$ and to ensure that the values in these adjacent cells
have the strongest control over $\phi_{f}$, rows associated with
these values in the least squares fit matrix are weighted a factor
of 1000 relative to the other rows \citep[following][]{Las02}. This
does not affect the order of accuracy. Mathematically, an arbitrarily
large value of the weight can be used to ensure that the fit goes
exactly through the upwind and downwind cell. However if a value too
large is used, the singular value problem becomes ill conditioned.
We do not use any of the other stabilisation procedures as described
by \citet{SWMD1x}. The value $\phi_{f}$ is then calculated as a
higher-order correction to first-order upwind: 
\begin{equation}
\phi_{f}=\phi_{\text{up}}+\sum_{c}w_{c}\phi_{c}\label{eq:multiDflux}
\end{equation}
where $w_{c}$ are the weights for each cell of the stencil calculated
from the least squares fit, with $w_{\text{up}}$ reduced by 1 to
make the fit a correction on upwind.

\subsubsection{Implicit Solution, Matrix Solvers and Tolerances \label{sub:implicitMethod}}

The trapezoidal implicit or Crank-Nicholson time-stepping leads to
a matrix equation which needs to be solved to find all the $\phi$s
at the next time-step. In order to ensure that the matrix is diagonally
dominant for arbitrary time-steps, the cubic interpolation applied
is a deferred correction on first-order upwind so that only the coefficient
corresponding to the upwind cells are included in the matrix. This
means that more than one implicit solves are needed per time-step
so that the higher-order terms are solved to be second-order accurate
in time. If the Courant number is less than or close to one, we use
two implicit solves per time-step. Consequently, assuming that the
velocity field and mesh are constant in time, the time-stepping scheme
is defined as:
\begin{eqnarray}
\frac{\phi^{\prime}-\phi^{n}}{\Delta t} & = & -\frac{1}{2V}\sum_{f\in c}\left(\phi_{\text{up}}^{n}+\phi_{\text{up}}^{\prime}+2\sum_{c}w_{c}\phi_{c}^{n}\right)\mathbf{u}_{f}\cdot\mathbf{S}_{f}\label{eq:CN_1}\\
\frac{\phi^{n+1}-\phi^{n}}{\Delta t} & = & -\frac{1}{2V}\sum_{f\in c}\left(\phi_{\text{up}}^{n}+\phi_{\text{up}}^{n+1}+\sum_{c}w_{c}\phi_{c}^{n}+\sum_{c}w_{c}\phi_{c}^{\prime}\right)\mathbf{u}_{f}\cdot\mathbf{S}_{f}.\label{eq:CN_2}
\end{eqnarray}
For larger Courant numbers we use four implicit solves per time-step
although sensitivity to this choice has not been investigated. 

If this implicit scheme is applied on a logically rectangular, two-dimensional
mesh with horizontal and vertical Courant numbers $c_{x}$ and $c_{z}$,
then the diagonal coefficients of the matrix would be $1+c_{x}/2+c_{z}/2$
and, assuming two upwind directions, there would be exactly two off-diagonal
elements, $-c_{x}/2$ and $-c_{z}/2$. Consequently the matrix is
very sparse, asymmetric and diagonally dominant for all time-steps.
It is solved using the OpenFOAM bi-conjugate gradient solver using
DILU pre-conditioning to a tolerance of $10^{-8}$ every iteration.
Sensitivity to the solver or solver tolerance have not been investigated.
Information about the number of solver iterations for different test
cases in given in section \ref{sub:cost}.

\section{Results of Test Cases in Planar Geometry\label{sec:results}}

In order to make test cases as simple as possible without the need
for incorporating spherical geometry, multi-panel meshes or non-rectangular
cells, our computational domain consists of a periodic two-dimensional
plane with deformations in the co-ordinate system (or mesh) mimicking
the kind of distortions which are produced by a cubed-sphere mesh.
We also use a two-dimensional vertical slice test case over orography
using terrain following co-ordinates (or terrain following meshes).

A range of test cases are undertaken on uniform and distorted meshes
using the dimensionally split scheme and the multi-dimensional scheme
in order to evaluate the influences of mesh distortions, the validity
of using a dimensionally split scheme on a distorted mesh and the
schemes' accuracy and stability for long time-steps.

\subsection{Solid Body Rotation\label{sub:solidRotate}}

The solid body rotation test case of \citet{LLM96} is used to compare
the accuracy of the dimensionally split and multi-dimensional schemes
on orthogonal and non-orthogonal meshes. We define this test case
on a domain that is $10^{4}\times10^{4}\text{ m}^{2}$. The velocity
is defined by numerically differentiating the streamfunction which
can be defined at mesh vertices, $\mathbf{x}$, by: 
\begin{equation}
\psi(\mathbf{x},t)=A|\mathbf{x}-\mathbf{x}_{c}|^{2}
\end{equation}
where $\mathbf{x}_{c}$ is the centre of the domain, and $A=5\pi/3000\text{ s}^{-1}$
so that the angular velocity is $2A$. The initial tracer takes a
Gaussian distribution in order to ensure that all advection schemes
achieve their theoretical order of accuracy:
\begin{equation}
\phi(\mathbf{x})=\exp\left(-\frac{1}{2}\frac{|\mathbf{x}-\mathbf{x}_{\phi}|^{2}}{r_{\phi}^{2}}\right)\label{eq:7}
\end{equation}
where $\mathbf{x}_{\phi}=\mathbf{x}_{c}+r_{c\phi}\mathbf{j}$ is the
initial centre of the tracer distribution, $r_{c\phi}=2500\text{ m}$,
$r_{\phi}=500\text{ m}$ and $\mathbf{j}$ is the unit vector in the
$y$ direction. The analytic solution has the same tracer distribution
but with the centre of the tracer at: 
\begin{equation}
\mathbf{x}_{\phi}=\mathbf{x}_{c}+r_{c\phi}\left(\begin{array}{c}
\cos(\pi/2+2At)\\
\sin(\pi/2+2At)
\end{array}\right)
\end{equation}
and the Gaussian rotates anti-clockwise exactly one revolution in
$600\text{ s}$.

\begin{figure}
\begin{centering}
\includegraphics[width=0.49\textwidth,height=0.49\textwidth]{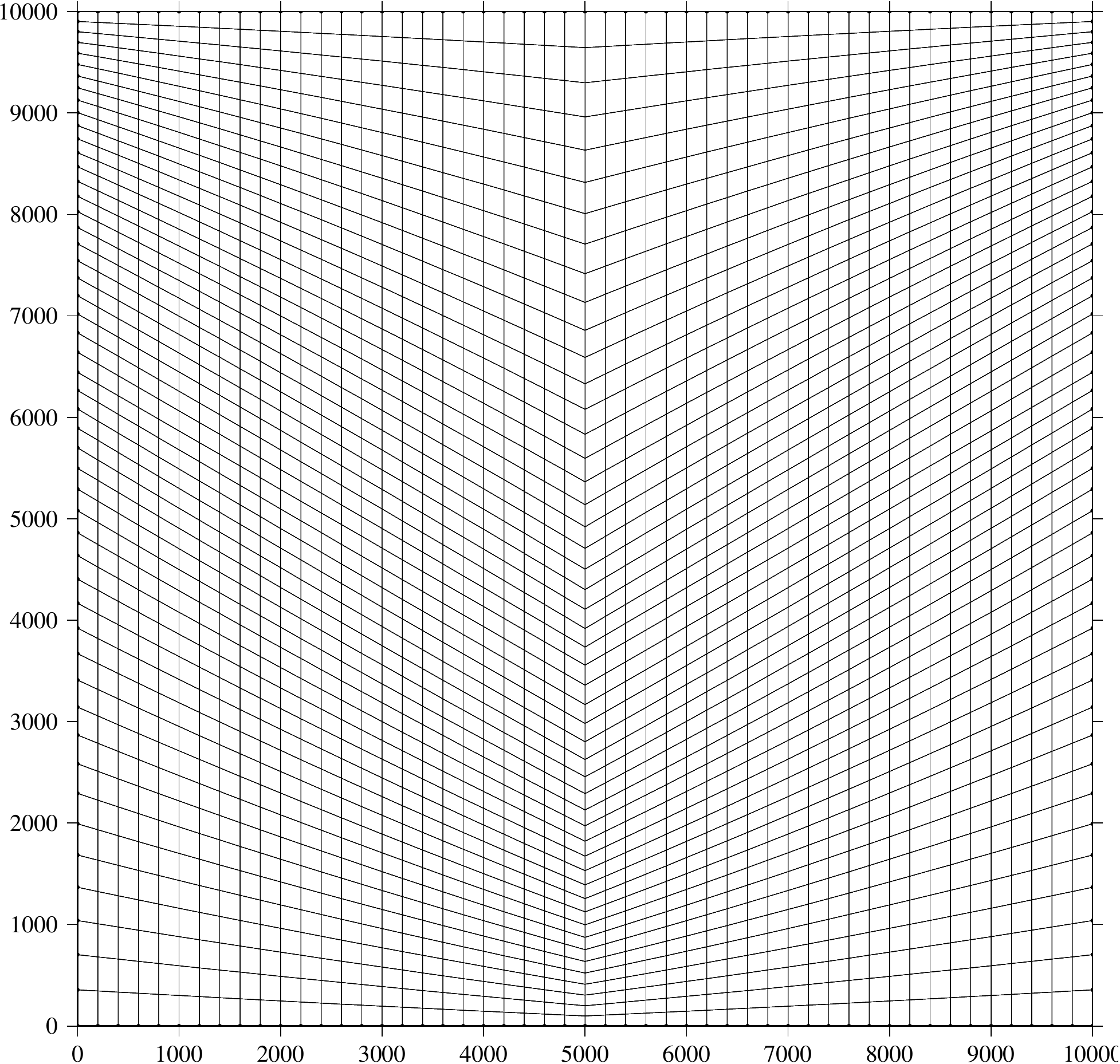}
\par\end{centering}

\caption{Two non-orthogonal mesh on a plane with $50\times50$ mesh points
and angles up to $120^{\text{o}}$ designed to have similar angles
to a cubed-sphere mesh. \label{fig:cubedGrid}}
\end{figure}

The solid body rotation test case is performed on a uniform, orthogonal
mesh and on a non-orthogonal mesh on a plane with non-orthogonality
similar to that of a cubed-sphere mesh. For the dimensionally split
scheme, non-orthogonality is achieved using the co-ordinate transform:
\begin{equation}
\begin{array}{ccc}
X=x &  & Y=\begin{cases}
y_{m}\left(1+\frac{y-f}{2y_{m}-f}\right) & \text{ for }y\ge f\\
y_{m}\left(1+\frac{y-f}{f}\right) & \text{ for }y<f
\end{cases}\end{array}
\end{equation}
where $y_{m}=5000\text{ m}$ and $f$ is the equation for a $y$ position
of uniform $Y$ half way up the domain. In order to create angles
of $120^{o}$ in the mesh, similar to a cubed-sphere, $f$ is given
by:
\begin{equation}
f=\begin{cases}
y_{m}\left(1+\frac{1}{2\sqrt{3}}\right)-\frac{x}{\sqrt{3}} & \text{ for }x\le x_{m}\\
y_{m}\left(1-\frac{1}{2\sqrt{3}}\right)+\frac{x-x_{m}}{\sqrt{3}} & \text{ for }x>x_{m}
\end{cases}.\label{eq:f_for_Vmesh}
\end{equation}
where $x_{m}=5000$ m. For a $50\times50$ mesh this gives the $x$
and $y$ co-ordinate locations as shown in figure \ref{fig:cubedGrid}.
The multi-dimensional scheme model uses Cartesian co-ordinates and
a distorted mesh rather than a non-orthogonal co-ordinate system on
a Cartesian mesh. However this does not affect the numerical results
assuming that the co-ordinate transforms are implemented in a consistent
way to the distorted mesh in Cartesian co-ordinates. 

For the dimensionally split scheme, bi-periodic boundary conditions
are applied. For the multi-dimensional scheme, it was more straightforward
to impose fixed value boundary conditions of $\phi=0$  where the
velocity is into the domain and zero normal gradient where the velocity
is out of the domain. However $\phi$ remains almost zero near the
boundaries so these boundary conditions do not affect the accuracy.

\begin{figure}
\includegraphics{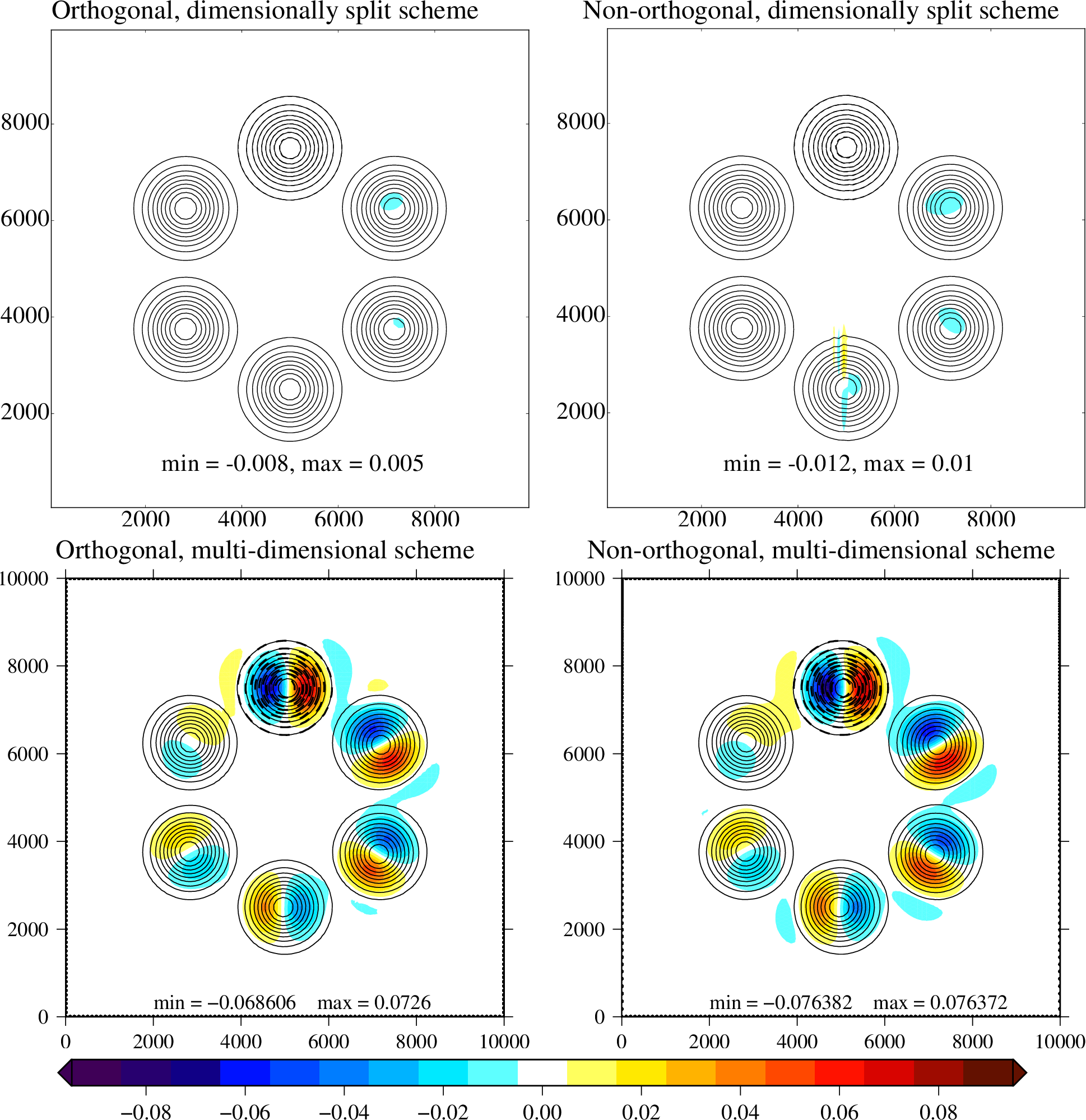}

\caption{Results of the solid body rotation test case on orthogonal and non-orthogonal meshes every 100 seconds for meshes of $100\times100$ cells and $\Delta t=1\text{ s}$ giving a maximum Courant number of around 1. The contours show the tracer value from 0.1 every 0.1 to 0.9 and the colours show the errors. The dashed contour shows the initial conditions. The min and max text in each plot gives the minimum and maximum errors. \label{fig:solidBody}}
\end{figure}

Results of this test case on the orthogonal and non-orthogonal meshes
of $100\times100$ cells with $\Delta t=1\text{ s}$ are shown in
figure \ref{fig:solidBody} for both advection schemes (which gives
a maximum Courant number close to one). The contours show the tracer
value every 100 seconds and the colours show the errors summed every
100 seconds. The dimensionally split scheme outperforms the multi-dimensional
scheme on both meshes due to the higher-order accuracy of the split
scheme. The dimensionally split scheme introduces a small error at
300 seconds where the tracer goes through the change in direction
of the mesh which would be ameliorated if we were using monotonicity
constraints. The second-order, multi-dimensional scheme shows phase
lag but errors are almost entirely insensitive the mesh distortions. 

\begin{figure}
\includegraphics{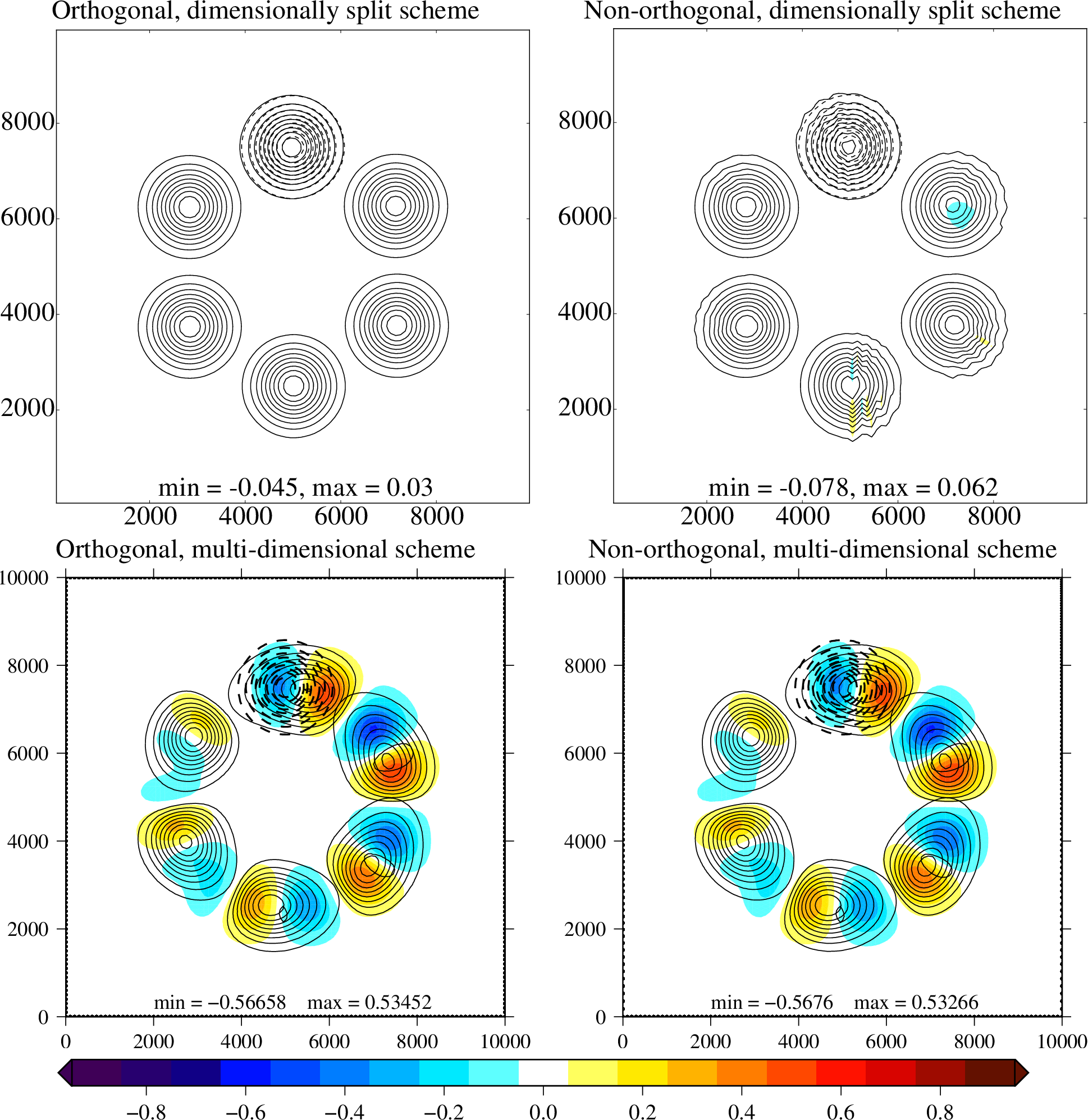}

\caption{Results of the solid body rotation test case on orthogonal and non-orthogonalmeshes every 100 seconds for meshes of $100\times100$ cells and $\Delta t=10\text{ s}$ giving a maximum Courant number of around 10. The contours show the tracer value from 0.1 every 0.1 to 0.9 and the colours show the errors. The dashed contour shows the initial conditions. The min and max text in each plot gives the minimum and maximum errors. \label{fig:solidBody-dt10}}
\end{figure}

The multi-dimensional and dimensionally split schemes take very different
approaches to handling large Courant numbers. The multi-dimensional
scheme uses implicit time-stepping whereas the dimensionally split
scheme uses a flux-form semi-Lagrangian approach, integrating over
a line of cells in order to calculate the flux across a face. Implicit
schemes are known to suffer from phase errors \citep[eg][]{DB12,LWW14}
for long time-steps whereas the accuracy of semi-Lagrangian is less
sensitive to time-step \citep{PS84}. Therefore we present results
of both schemes on orthogonal and non-orthogonal meshes for time-steps
10 times those used in figure \ref{fig:solidBody} ($\Delta t=10\text{ s}$)
giving maximum Courant numbers of around 10 using $100\times100$
cells in figure \ref{fig:solidBody-dt10}. The error of the multi-dimensional
scheme is again much larger than the dimensionally split scheme on
both meshes. The dimensionally split scheme is accurate at large Courant
numbers despite the first-order calculation of departure points. However,
the dimensionally split scheme introduces oscillations on the non-orthogonal
mesh, particularly where the mesh changes direction whereas for the
multi-dimensional scheme, errors are not strongly affected by the
non-orthogonality. Figure \ref{fig:solidBody-dt10} clearly shows
phase errors of the implicit time-stepping but, despite the large
Courant number, the well resolved part of the profile is propagating
at close to the correct speed. Dispersion analysis \citep{LWW14}
shows that high \emph{frequency} oscillations (which are poorly resolved
in time \emph{and} space) will be slowed dramatically but fast moving
features which are \emph{well resolved in space} will propagate at
a much more realistic speed, supporting the results shown in figure
\ref{fig:solidBody-dt10}. 

\begin{figure}
\begin{centering}
\includegraphics{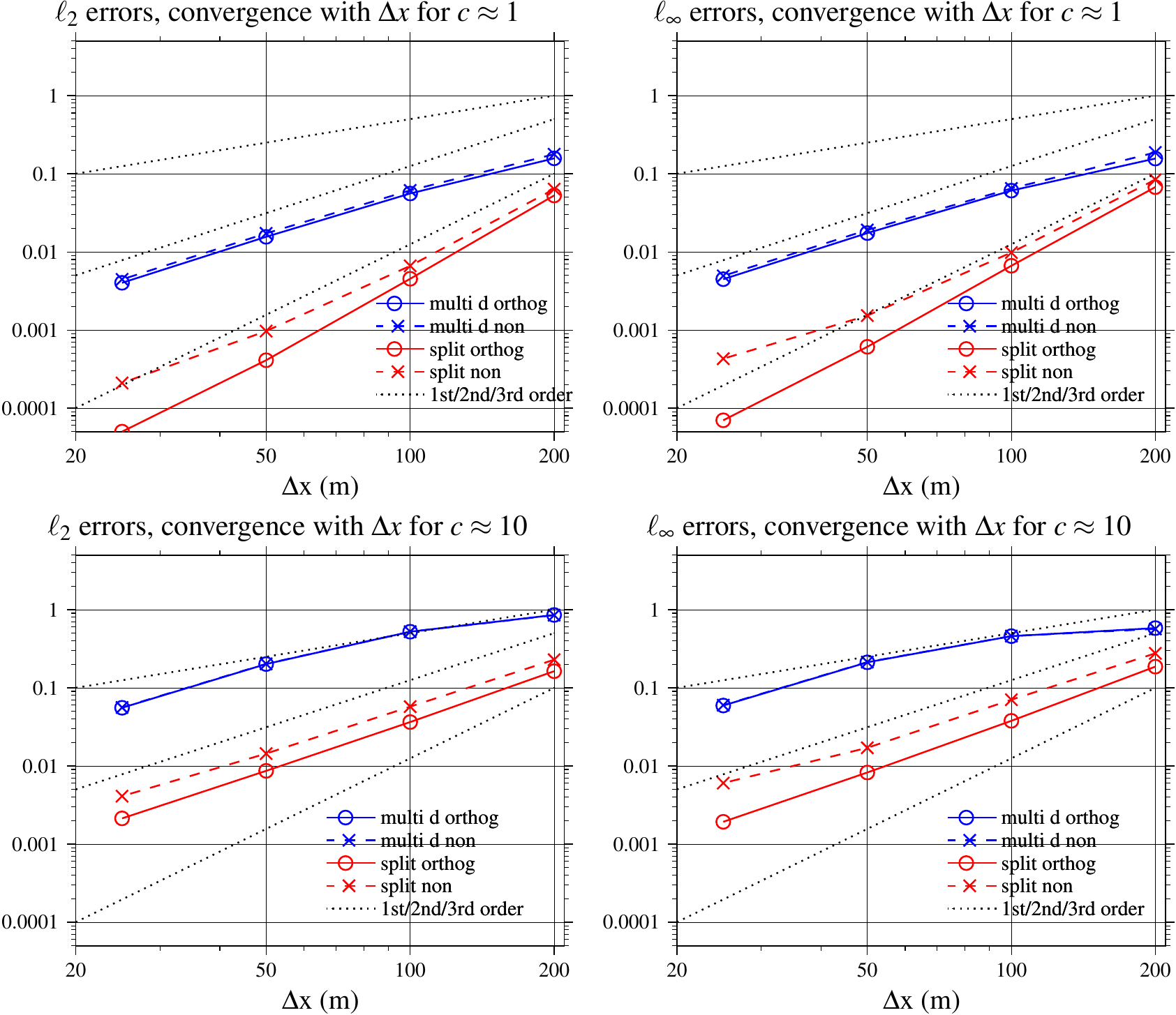}
\par\end{centering}

\caption{Convergence with spatial resolution of the $\ell_{2}$ and $\ell_{\infty}$ errors for the solid body rotation test case at $t=500\ \text{seconds}$ ($5/6$ of one revolution) on orthogonal and non-orthogonal meshes
using multi-dimensional and dimensionally split schemes. \label{fig:solidBodyConverge} }
\end{figure}

In order to compare convergence with resolution of the different numerical
methods, we use the $\ell_{2}$ and $\ell_{\infty}$ error norms defined
in the usual way:
\begin{eqnarray}
\ell_{2} & = & \sqrt{\int_{V}\left(\phi-\phi_{T}\right)^{2}dV\bigg/\int_{V}\phi_{T}^{2}\ dV}\\
\ell_{\infty} & = & \max|\phi-\phi_{T}|\big/\max|\phi_{T}|
\end{eqnarray}
where $\phi_{T}$ is the analytic solution and the integrations and
maxima are over the whole domain, with volume $V$. Figure \ref{fig:solidBodyConverge}
shows convergence with resolution of the $\ell_{2}$ and $\ell_{\infty}$
error measures for meshes of $50\times50$, $100\times100$, $200\times200$
and $400\times400$ cells with time-steps scaled in order to maintain
a maximum Courant number of 1 ($\Delta t=2,\ 1,\ 0.5,\ 0.25\text{ s}$)
or scaled to achieve a maximum Courant number of 10 ($\Delta t=20,\ 10,\ 5,\ 2.5\text{ s}$).
The error norms are calculated at $t=500\ \text{seconds}$, when the
tracer has made $5/6$ of one revolution in order to avoid error cancellation.
On both orthogonal and non-orthogonal meshes, the multi-dimensional
scheme has second order convergence once errors are low enough to
avoid error saturation (stable errors are bounded at around one).
With a large Courant number, both schemes are less accurate. For the
multi-dimensional scheme, this is due to phase errors of the implicit
time-stepping whereas for the multi-dimensional scheme the first-order
errors in calculating the departure point and trajectory will be emerging
and there could also be significant errors from the second-order COSMIC
splitting. On the orthogonal mesh, the dimensionally split scheme
has  third order convergence for the Courant number close to one and
second order for the larger Courant number \citep[consistent with the results of ][]{CW84,LLM96}.
On the non-orthogonal mesh, the dimensionally split scheme has second
order converges for both Courant numbers.

\begin{figure}
\begin{centering}
\includegraphics{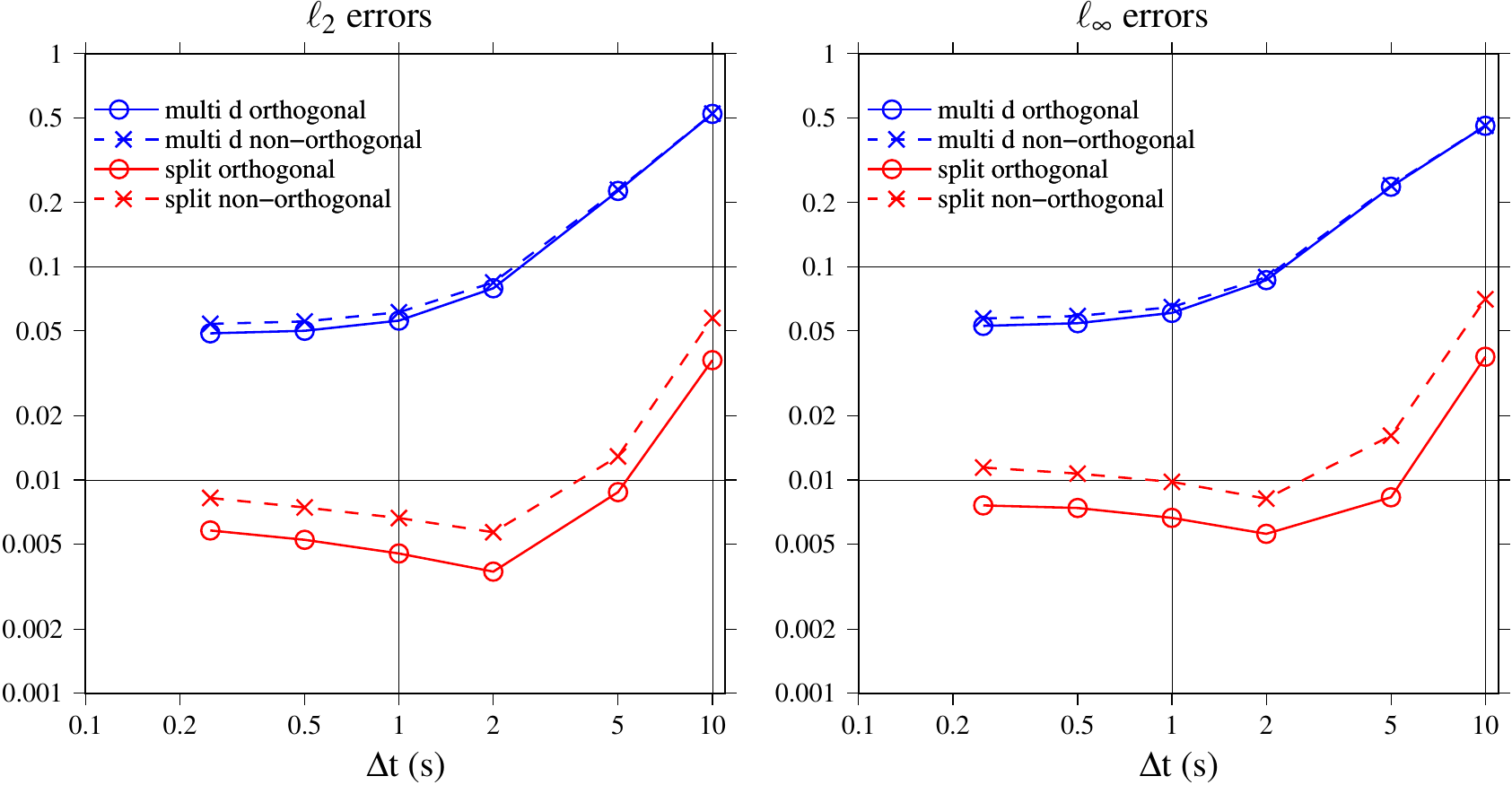}
\par\end{centering}

\caption{Convergence with temporal resolution of the $\ell_{2}$ and $\ell_{\infty}$errors for the solid body rotation test case at $t=500\ \text{seconds}$ ($5/6$ of one revolution) using orthogonal and non-orthogonal meshes of $100\times100$ cells using multi-dimensional and dimensionally split schemes. Maximum Courant numbers range from $c=0.25$ (for $\Delta t=0.25\ \text{s})$ to $ $$c=10$ (for $\Delta t=10\ \text{s})$. \label{fig:solidBodyConverge-dt} }
\end{figure}

The different time-stepping schemes of the two models also affect
accuracy. The $\ell_{2}$ and $\ell_{\infty}$ error measures as a
function of time-step for meshes of $100\times100$ cells are shown
in figure \ref{fig:solidBodyConverge-dt}. The dimensionally split
scheme, which uses flux-form semi-Lagrangian time-stepping, has errors
reducing as time-step increases, up to a Courant number of 2 ($\Delta t=2\ \text{s}$)
whereas the multi-dimensional scheme, which uses the method of lines
to treat space and time separately, always has error reducing as time-steps
reduces. The flux-form semi-Lagrangian technique discretises space
and time together and the error is not very sensitive to time-step.
However, the shorter the time-step, the more time-steps need to be
taken so errors can actually accumulate more by taking more time-steps.
(This is consistent with the order of accuracy of semi-Lagrangian
being $\Delta x^{p}/\Delta t$ for interpolation using polynomials
of degree $p$, as described by \citet{Dur10}.)

In summary, the dimensionally split scheme has excellent behaviour
at large and small Courant numbers on the orthogonal meshes with up
to third order convergence for small Courant numbers and the errors
increase and order of convergence decreases on non-orthogonal meshes.
In contrast, the multi-dimensional scheme is insensitive to the orthogonality,
converges with second-order and suffers from phase errors at large
Courant numbers.

\subsection{Horizontal Advection over Orography\label{sub:overOrog}}

Non-orthogonal meshes (or co-ordinate systems) are usually necessary
for representing orography which could be a challenge for dimensionally
split schemes. Horizontal-vertical split schemes are commonly used
in this context \citep[eg.][]{DEE+12,WGZ+13}. We present results
of the \citet{SLF+02} horizontal advection over orography test case
for a range of resolutions and Courant numbers for the dimensionally
split and multi-dimensional schemes.

All simulations use basic terrain following co-ordinates \citep[BTF,][]{GCS75}
in order to present a challenging test case that maximises the non-orthogonality.
The transformation is given by:
\begin{equation}
\begin{array}{ccc}
X=x &  & Z=H\frac{z-h(x)}{H-h(x)}\end{array}
\end{equation}
where $H$ is the domain height and $h$ is the terrain height. The
test case uses a domain of width 300\ km, height, $H=25\text{ km}$
and a mountain range defined by
\begin{equation}
h=\begin{cases}
h_{0}\cos^{2}\frac{\pi x}{\lambda}\cos^{2}\frac{\pi x}{2a} & \text{ for }|x|\le a\\
0 & \text{ otherwise}
\end{cases}
\end{equation}
with the maximum mountain height, $h_{0}=3\text{ km}$, half-width
$a=25\text{ km}$ and wavelength $\lambda=8\text{ km}$. These values
give a maximum terrain gradient of close to $45^{o}$. The wind is
given by a streamfunction which is defined at vertices so that the
wind field is discretely divergence free. The streamfunction at vertices
is calculated analytically from the wind profile:
\begin{equation}
u(z)=u_{0}\begin{cases}
1 & \text{ for }z_{2}\le z\\
\sin^{2}\left(\frac{\pi}{2}\frac{z-z_{1}}{z_{2}-z_{1}}\right) & \text{ for }z_{1}<z\le z_{2}\\
0 & \text{ for }z<z_{1}
\end{cases}
\end{equation}
with $u_{0}=10\text{ ms}^{-1}$, $z_{1}=4\text{ km}$ and $z_{2}=5\text{ km}$.
The initial tracer position is given by:
\begin{equation}
\phi=\begin{cases}
\cos^{2}\frac{\pi r}{2} & \text{ for }r\le1\\
0 & \text{ otherwise}
\end{cases}\text{ with }r=\sqrt{\left(\frac{x-x_{0}}{A_{x}}\right)^{2}+\left(\frac{z-z_{0}}{A_{z}}\right)^{2}}
\end{equation}
with initial tracer centre, $(x_{0},z_{0})=(-50\text{ km},\ 9\text{ km})$
and halfwidths $A_{x}=25\text{ km}$, $A_{z}=9\text{ km}$. At time
$t=5000\text{ s}$ the tracer is above the mountain and the simulation
finishes at $t=10,000\text{ s}$ by which time the analytic solution
is centred at $(50\text{km },8\text{km})$.

\begin{figure}
\includegraphics{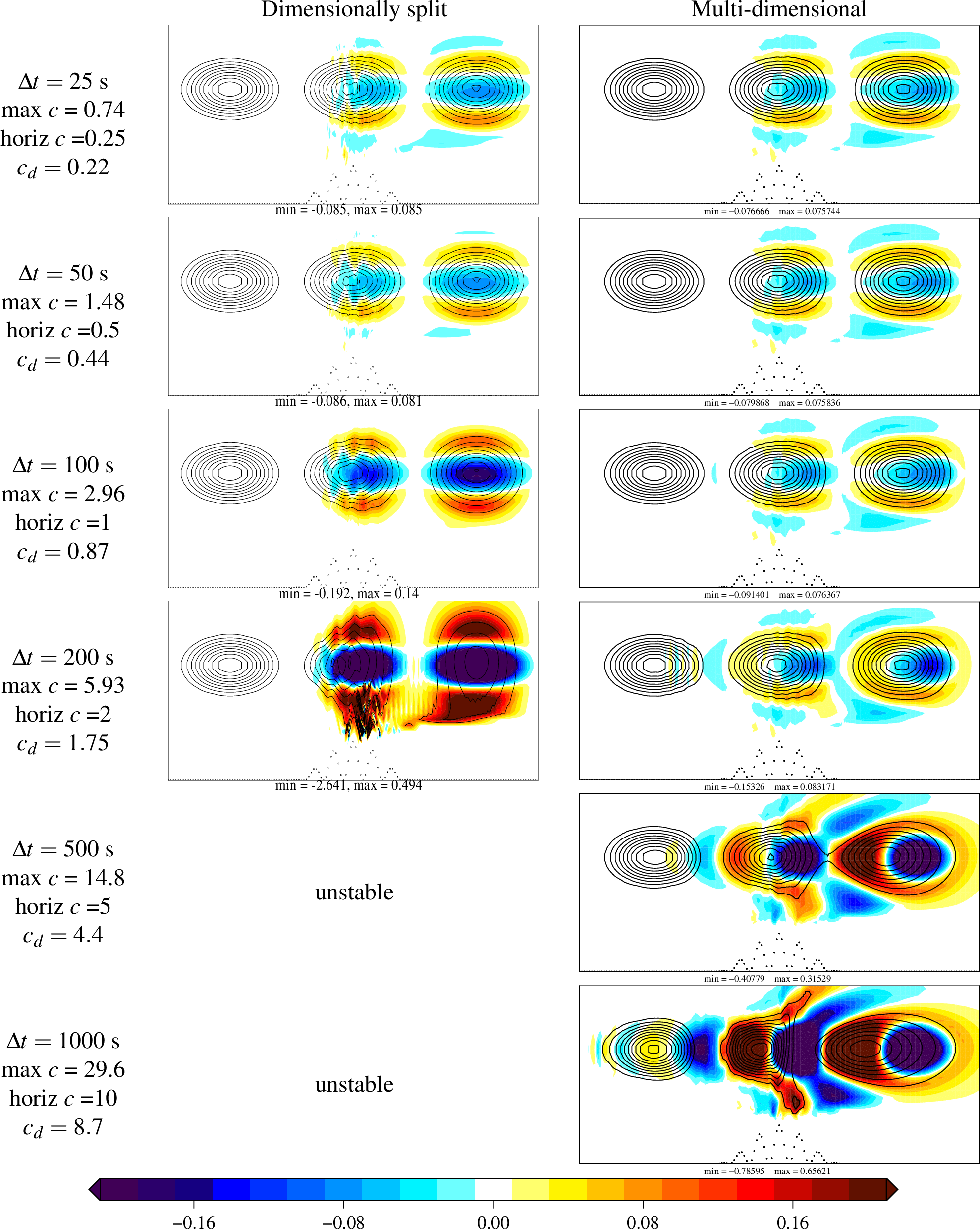}

\caption{Tracer contours after 0, 5000 and 10000s for the horizontal advection over orography and errors in colour. Different plots use different time-steps. Spatial resolution is $\Delta x=1\ $km, $\Delta z=500\ $m. Maximum Courant numbers, horizontal Courant number, $c$, and maximum deformational Courant number, $c_{d}$ are also given. \label{fig:horizontalAdvectionWithDt}}
\end{figure}

The tracer advection over orography is shown in figure \ref{fig:horizontalAdvectionWithDt}
for the split and multi-dimensional schemes at a resolution of $\Delta x=1\text{ km}$,
$\Delta z=500\text{ m}$ and for a range of Courant numbers. The horizontal
Courant number is defined as $u_{0}\Delta t/\Delta x$ and ranges
from 0.25 to 10. The maximum Courant number is the maximum of the
multi-dimensional Courant number which is defined for cell with faces
$f$ as:
\begin{equation}
c=\frac{1}{2V}\sum_{f}|\mathbf{u}_{f}\cdot\mathbf{S}_{f}|\Delta t
\end{equation}
(see section \ref{sub:MMmolscheme} for definitions of variables)
and ranges from 0.74 to 29.6. The time-step restriction for the split
scheme is based on the deformational Courant number, $c_{d}\le1$,
(eqn \ref{eq:cd}). The maximum deformational Courant number is also
given in figure \ref{fig:horizontalAdvectionWithDt}. The contours
in figure \ref{fig:horizontalAdvectionWithDt} show the tracer values
at 0, $5000\text{ s}$ and $10,000\text{ s}$ after initialisation
and the colours show the errors from the analytic solution. For horizontal
Courant numbers less than one (maximum Courant number up to 3), both
schemes give accurate results with the dimensionally split scheme
tending to give oscillations and the multi-dimensional scheme producing
more diffusion. For larger Courant numbers, when the deformational
Courant number is greater than one, the split scheme is unstable while
the multi-dimensional scheme produces large phase errors due to the
errors associated with implicit time-stepping. The term responsible
for the large deformational Courant number is $\partial u/\partial z$
where the velocity shears from $u_{0}=10$m/s at $z=z_{2}$ to zero
at $z=z_{1}$. 

\begin{figure}
\noindent \begin{centering}
\includegraphics{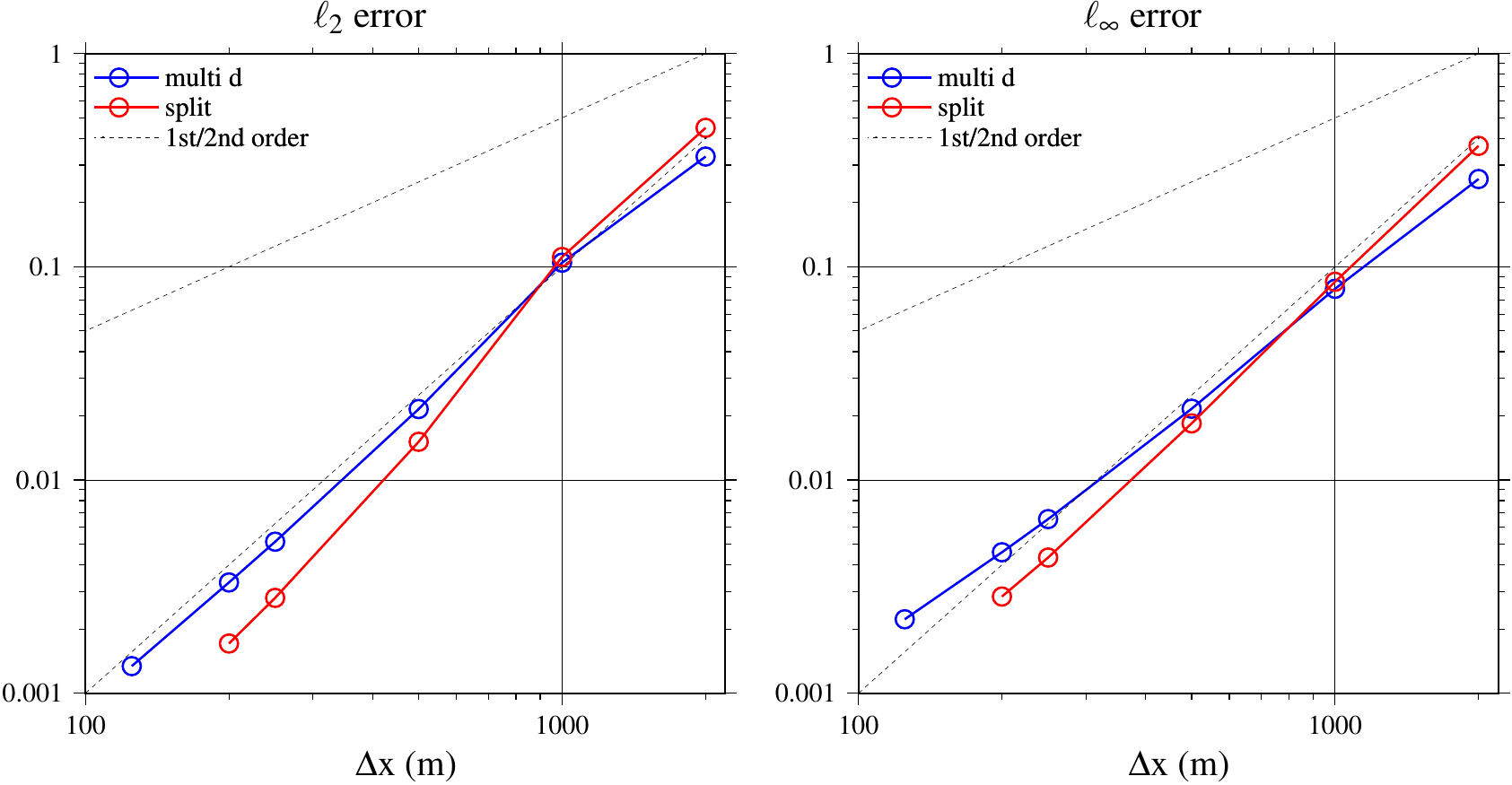}
\par\end{centering}

\caption{Convergence of the $\ell_{2}$ and $\ell_{\infty}$ error metricsfor horizontal advection over a 3000m and a 6000m mountains after 10,000s for the dimensionally split (PPM) scheme and multi-dimensional MOL scheme, all using a horizontal Courant number of 0.25. Lines also show theoretical first- and second-order convergence are also shown.
\label{fig:horizontalAdvectConverge}}
\end{figure}

We examine the convergence with resolution in figure \ref{fig:horizontalAdvectConverge}
which shows the $\ell_{2}$ and $\ell_{\infty}$ error norms as a
function of $\Delta x$. These simulations all use a maximum Courant
number less than one, a horizontal Courant number of 0.25, a maximum
deformational Courant number of about 0.2 and fixed ratios of $\Delta x$,
$\Delta z$ and $\Delta t$. Both schemes give similar accuracy with
the dimensionally split scheme having faster convergence with resolution.

In summary, the dimensionally split scheme has good accuracy over
orography for modest Courant numbers but larger errors for larger
Courant numbers and is unstable when the deformational Courant number
is greater than one whereas the multi-dimensional scheme with implicit
time-stepping is stable for all Courant numbers. The multi-dimensional
scheme second-order convergent and the dimensionally split scheme
converges faster.

\subsection{Deformational Flow \label{sub:deformFlow}}

In deformational flow, there is no analytical solution and therefore
we follow the approach of \citet{NL10,LSPT12} and define an evolving
velocity field that reverses direction half way through the simulation,
taking the tracer back to the initial conditions. Error norms can
then be calculated by comparing the final and initial tracer fields.
The deformational velocity field is added to a fixed solid body rotation
that does not reverse so as to avoid error cancellation between the
forwards and backwards periods (solid body here meaning horizontal
wind on a periodic domain). 

We define a Cartesian version of the deformational, non-divergent
test case of \citet{LSPT12} with a domain between $-\pi$ and $\pi$
in the $x$ direction ($L_{x}=2\pi$) and between $-\frac{\pi}{2}$
and $\frac{\pi}{2}$ in the $y$ direction ($L_{y}=\pi$) with periodic
boundary conditions in the $x$ direction and zero gradient, zero
flow boundary conditions in the $y$ direction. The stream function
adapted to Cartesian co-ordinate is:
\begin{eqnarray}
\psi\left(x,y,t\right) & = & \frac{\hat{\psi}}{T}\left(\frac{L_{x}}{2\pi}\right)^{2}\sin^{2}\left(2\pi\left(\frac{x}{L_{x}}-\frac{t}{T}\right)\right)\cos^{2}\left(\pi\frac{y}{L_{y}}\right)\cos\left(\pi\frac{t}{T}\right)-L_{x}\frac{y}{T}
\end{eqnarray}
where $\hat{\psi}=10$ and $T=5$ is the time for one complete revolution
of the periodic domain for the solid body rotation part of the flow.
In order to test the order of convergence of the advection schemes,
we use the infinitely smooth Gaussian distribution for the initial
tracer concentration:
\begin{equation}
\phi=0.95\exp-\frac{|\mathbf{x}-\mathbf{x}_{0}|^{2}}{A}+0.95\exp-\frac{|\mathbf{x}-\mathbf{x}_{1}|^{2}}{A}
\end{equation}
where $\mathbf{x}^{T}=(x,y)$, $\mathbf{x}_{0}^{T}=(\frac{5}{12}L_{x},0)$,
$\mathbf{x}_{1}^{T}=(\frac{7}{12}L_{x},0)$ and $A=\frac{1}{5}$.

\begin{figure}
\includegraphics{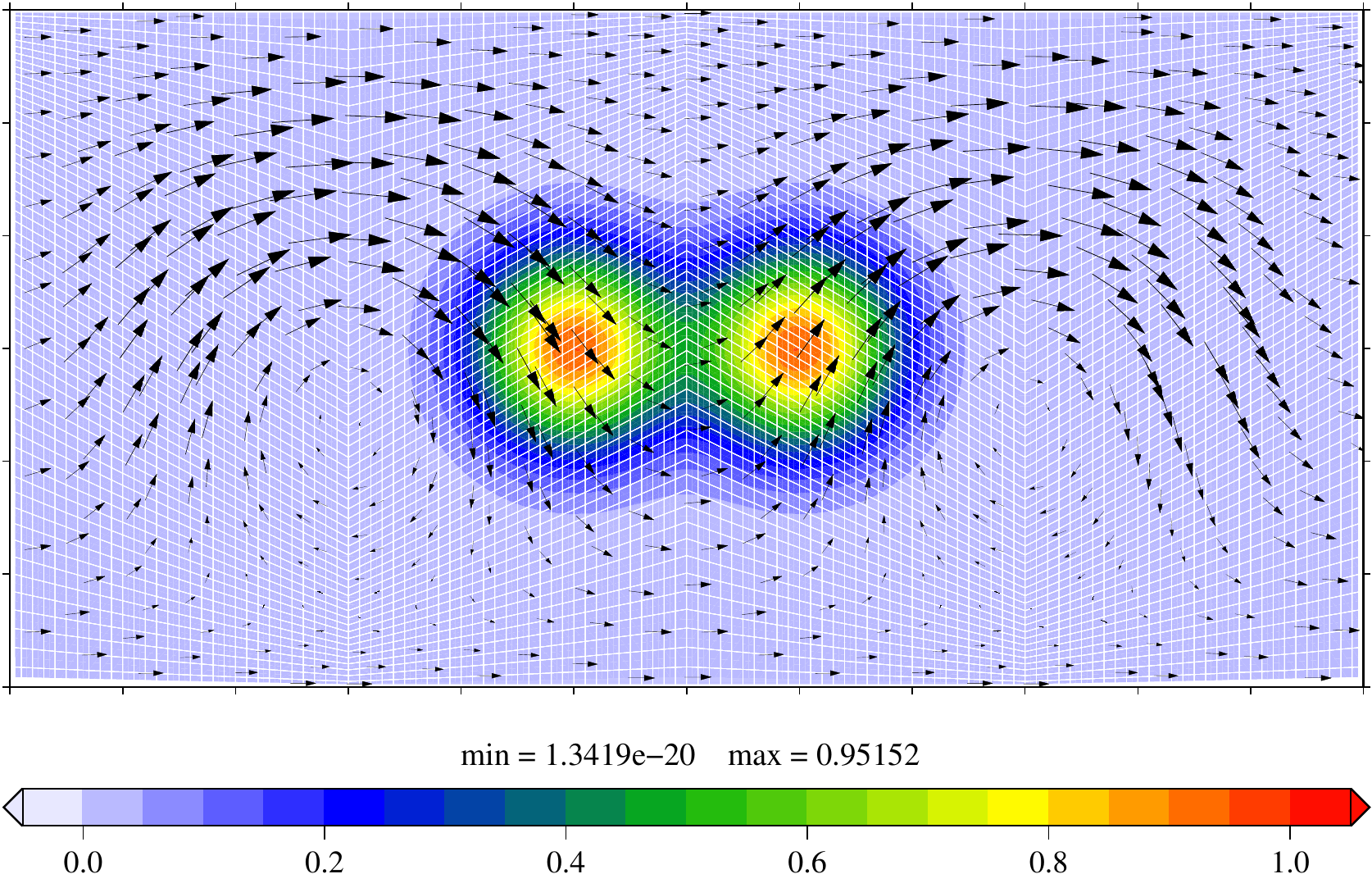}

\caption{The non-orthognal mesh/co-ordinate system for the deformational flow with $120\times60$ cells and the initial tracer conditions.\label{fig:deformGrid}}
\end{figure}

A mesh with distortion similar to a cubed-sphere is defined by the
co-ordinate transform:
\begin{equation}
\begin{array}{ccc}
X=x &  & Y=\begin{cases}
L_{y}\frac{y-f}{L_{y}-2f} & \text{ for }y\ge f\\
L_{y}\frac{y-f}{L_{y}+2f} & \text{ for }y<f
\end{cases}\end{array}
\end{equation}
where $f$ is:
\[
f=\begin{cases}
\frac{1}{\sqrt{3}}\left(\frac{\pi}{4}-|x|\right) & \text{ for }|x|\le\frac{\pi}{2}\\
\frac{1}{\sqrt{3}}\left(|x|-\frac{3\pi}{4}\right) & \text{ for }|x|>\frac{\pi}{2}
\end{cases}.
\]
The initial conditions and a $120\times60$ mesh with distortion defined
by this co-ordinate transform is shown in figure \ref{fig:deformGrid}. 

\begin{figure}
\includegraphics{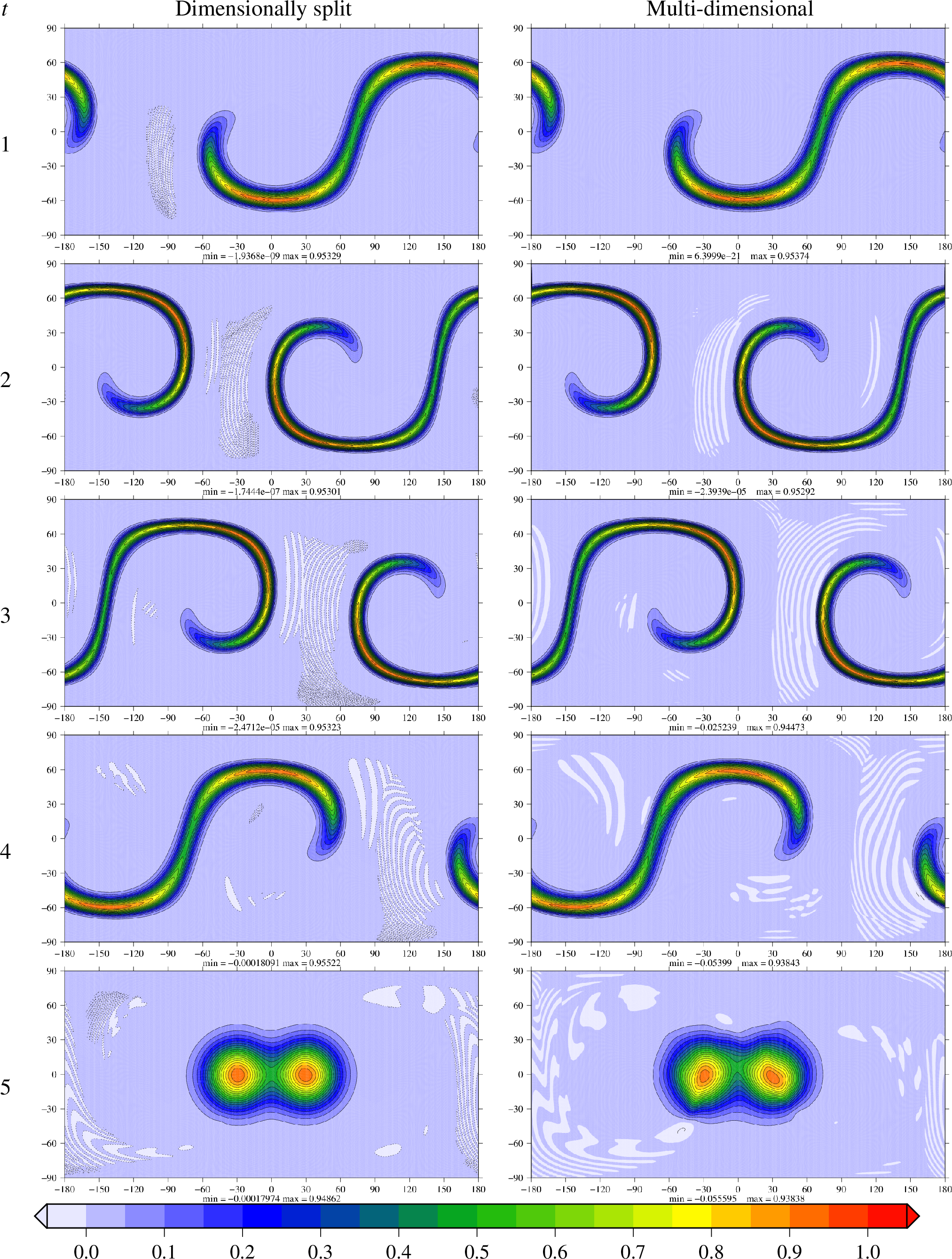}

\caption{Results of the deformational flow test case after 1, 2, 3, 4 and 5 time units on the non-orthogonal meshes of 480$\times$240 cells using the multi-dimensional and split schemes. \label{fig:deformResults}}
\end{figure}

The tracer concentrations after 1, 2, 3, 4 and 5 time units are shown
in figure \ref{fig:deformResults} using the dimensionally split and
multi-dimensional schemes on a non-orthogonal mesh of $480\times240$
cells, a time-step of 0.0025 units (ie 2000 time-steps to reach 5
time units) which gives a maximum Courant number of 1.03. The tracer
is stretched out, wound up, advected around and then wound back into
its original position with some numerical errors. Both advection schemes
preserve fine filaments but suffer from some dispersion errors which
generate small oscillations around zero behind sharp gradients in
the direction of the flow since neither scheme is monotonic or positive
preserving. The dimensionally split scheme returns the tracer to a
more accurate final solution.

\begin{table}
\begin{tabular}{|c||c|c|c|c|c|c|}
\hline 
Resolution & $60\times30$ & $120\times60$ & $240\times120$ & $480\times240$ & $960\times480$ & $1920\times960$\tabularnewline
\hline 
$\Delta t$ $\downarrow$ & \multicolumn{6}{c|}{Maximum Courant numbers {\footnotesize{}(deformational)} $\downarrow$}\tabularnewline
\hline 
\hline 
0.2 & 10.3 {\footnotesize{}(0.639)} &  &  &  &  & \tabularnewline
\hline 
0.1 &  & 10.3{\footnotesize{} (0.312)} &  &  &  & \tabularnewline
\hline 
0.05 &  & 5.17{\footnotesize{} (0.154)} & 10.3 {\footnotesize{}(0.154)} &  &  & \tabularnewline
\hline 
0.025 &  &  &  & 10.3{\footnotesize{} (0.077)} &  & \tabularnewline
\hline 
0.02 &  & 2.06{\footnotesize{} (0.061)} &  &  &  & \tabularnewline
\hline 
0.0125 & 1.03 &  &  &  & 10.3 & \tabularnewline
\hline 
0.01 &  & 1.03 {\footnotesize{}(0.030)} &  &  &  & \tabularnewline
\hline 
0.00625 &  &  &  &  &  & 10.3\tabularnewline
\hline 
0.005 &  & 0.517 {\footnotesize{}(0.015)} & 1.03{\footnotesize{} (0.015)} &  &  & \tabularnewline
\hline 
0.0025 &  & 0.254 {\footnotesize{}(0.008)} &  & 1.03{\footnotesize{} (0.008)} &  & \tabularnewline
\hline 
0.00125 &  &  &  &  & 1.03{\footnotesize{} (0.004)} & \tabularnewline
\hline 
0.000625 &  &  &  &  &  & 1.03\tabularnewline
\hline 
\end{tabular}

\caption{Resolutions, time-steps and Courant numbers for the non-orthogonal
meshes for the deformational flow test case. Maximum Courant numbers
are shown inside the table for simulations which are presented in
graphs. The maximum Courant numbers on the orthogonal meshes are 70\%
of those on the non-orthogonal meshes. The maximum deformational Courant
numbers are in brackets for some cases \label{tab:res_Courant}}
\end{table}

\begin{figure}
\includegraphics{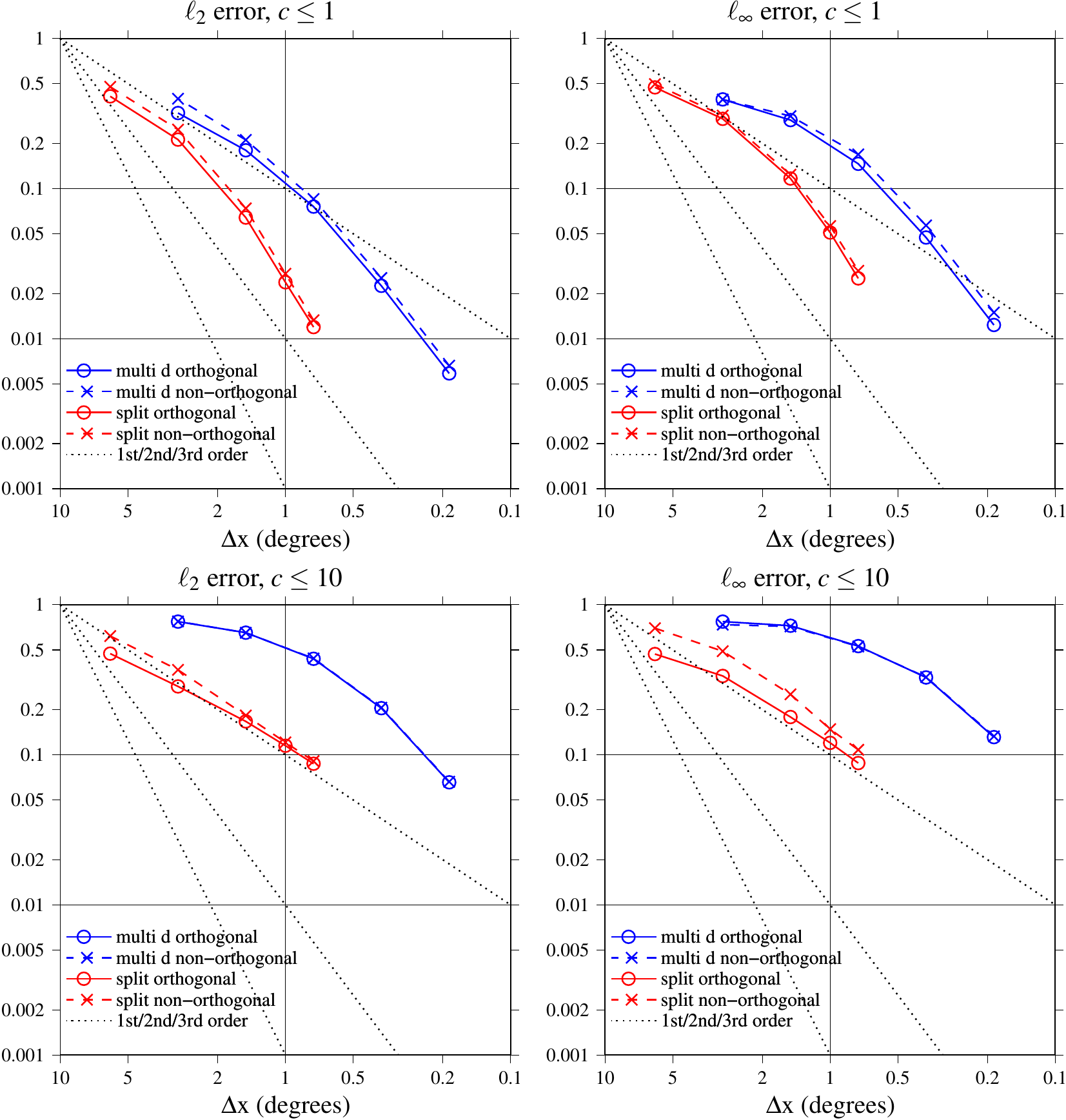}

\caption{Convergence with resolution of the $\ell_{2}$ and $\ell_{\infty}$ errors for the deformational flow on orthogonal and non-orthogonal meshes using multi-dimensional and dimensionally split schemes. \label{fig:deformationalConverge} }
\end{figure}

Sensitivity to orthogonality, resolution and time-step are explored
using a range of simulations with maximum Courant numbers shown in
table \ref{tab:res_Courant}. The Courant numbers on the uniform,
orthogonal meshes (with no mesh distortions) are 70\% of those on
the non-orthogonal mesh since the non-orthogonal mesh has clustering
of mesh points. The deformational Courant number is less than one
for all simulations. The convergence with resolution of the $\ell_{2}$
and $\ell_{\infty}$ error norms at the final time are shown in figure
\ref{fig:deformationalConverge} for both advection schemes using
a maximum Courant number close to one and for a maximum Courant number
close to ten. We will first consider the behaviour of the schemes
at modest Courant number (maximum close to one). All of the schemes
give first-order convergence with resolution at coarse resolutions
due to error saturation (if the scheme is stable, errors are bounded
above by about one). At higher resolution, the split scheme converges
with nearly third-order for both meshes whereas the multi-dimensional
scheme approaches second-order. For large Courant numbers (maximum
close to 10), the split scheme converges with first-order for both
mesh types due to the crude estimation of the departure point (section
\ref{sub:PPM}) whereas the multi-dimensional scheme is much less
accurate but approaches second-order  at high resolution. The dimensionally
split scheme is sensitive to mesh distortions only for large Courant
numbers. 

\begin{figure}
\includegraphics{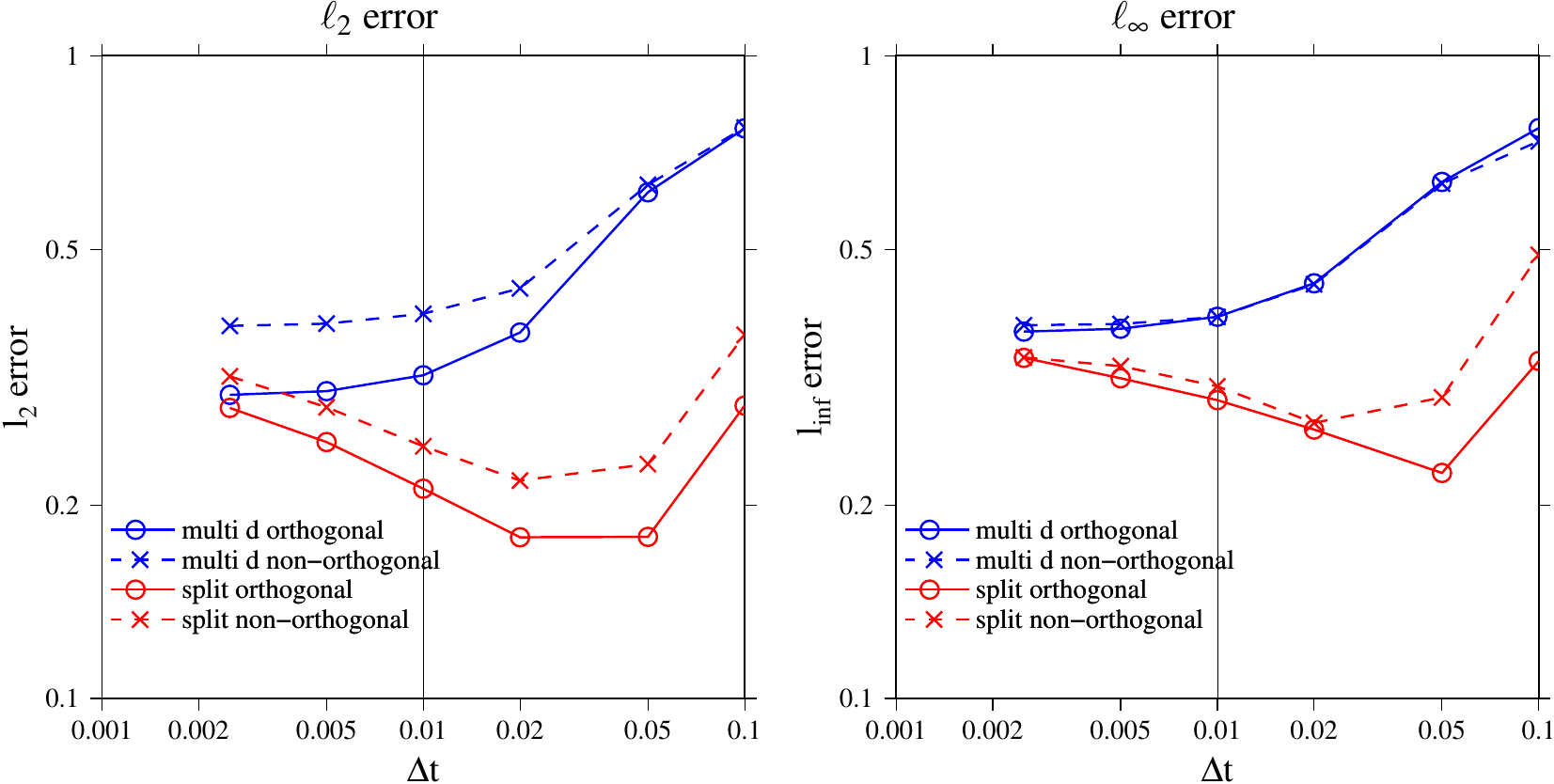}

\caption{Variation of $\ell_{2}$ and $\ell_{\infty}$ errors with time-step for the deformational flow on orthogonal and non-orthogonal meshes of $120\times60$ cells using multi-dimensional and dimensionally split schemes. Different time-steps give maximum Courant numbers of 0.25, 0.5, 1, 2, 5 and 10. \label{fig:deformational_dt} }
\end{figure}

We inspect sensitivity to time-step of both schemes on orthogonal
and non-orthogonal meshes of $120\times60$ cells in figure \ref{fig:deformational_dt}
using the time-steps shown in table \ref{tab:res_Courant} giving
maximum Courant numbers ranging from 0.25 to 10 (deformational Courant
numbers always less than 0.32). This demonstrates potentially useful
properties of the flux-form semi-Lagrangian time-stepping used by
the split scheme. Ignoring errors in calculating the departure point,
semi-Lagrangian schemes have errors proportional to $\Delta t^{-1}$
\citep{Dur10} which explains the reduction in error as time-step
increases for the split scheme. Once the maximum Courant number reaches
1 or 2 ($\Delta t=0.01$ or 0.02), errors of the split scheme do grow
with the time-step, due to the deformational nature of the flow and
errors in calculating the departure points. In contrast, using the
multi-dimensional scheme which uses method of lines time-stepping,
errors always increase as the time-step is increased. Comparisons
between schemes in figures \ref{fig:deformResults} and \ref{fig:deformationalConverge}
used maximum Courant numbers of 1 and 10 which showed that the split
scheme on the orthogonal mesh gave better accuracy than the multi-dimensional
scheme. However, figure \ref{fig:deformational_dt} shows that this
advantage disappears at lower Courant numbers since the multi-dimensional
scheme (method of lines) gets more accurate with lower Courant numbers
whereas the split scheme (semi-Lagrangian) gets less accurate since
more time-steps have to be taken. Both schemes are stable for all
time-steps considered.

In summary, the dimensionally split scheme has good accuracy for deformational
flow independent of mesh orthogonality. The multi-dimensional scheme
is competitive at small Courant numbers but the semi-Lagrangian nature
of the split scheme means that errors are very low for Courant numbers
close to one.

\subsection{Computational Cost\label{sub:cost}}

We cannot compare CPU time, wall clock time or parallel efficiency
of the two advection schemes because the multi-dimensional scheme
is written in C++ using OpenFOAM and the split scheme is written in
Python, both codes have been run on different hardware and the split
advection scheme code has not been parallelised. Instead we consider
the number of multiply operations performed per cell per time-step
for calculating the fluxes of each tracer. We appreciate that this
is not a good predictor of efficiency as it does not consider memory
read and write requirements or cache coherency. However all data that
is multiplied has to be fetched from memory and so, assuming that
all data can be arranged optimally in memory to enable fewest cache
misses, the number of multiplies should be related to the wall clock
time. Neither advection scheme does significantly more work per memory
fetch than the other. 

For both schemes we consider the number of multiplies needed to calculate
the flux of $\phi$ at each face. That is equation \ref{eq:PPMflux}
for the dimensionally scheme and equation (\ref{eq:multiDflux}) for
the multi-dimensional scheme. We do not consider the computational
cost of updating the cell averages from the fluxes using Gauss's theorem
(eqns (\ref{PPMgauss}) and (\ref{eq:multidGauss})) as these are
the same for each scheme. For the multi-dimensional scheme, we do
include an estimate of the amount of work done by the linear equation
solver but we do not consider the scalability of this solver.

\subsubsection{Dimensionally Split Scheme}

The computational cost of PPM with COSMIC splitting is not strongly
time-step dependent due to the swept area approach of flux-form semi-Lagrangian.
One-dimensional PPM uses 4 cells for interpolation to find face values
used by the reconstruction. Assuming as much as possible is pre-computed,
this interpolation uses 3 multiplies on a non-uniform mesh. The reconstruction
then uses 6 multiplies to calculate the flux on each face. Only one
additional memory access and one additional  multiply are needed per
cell for Courant numbers greater than one (eqn \ref{eq:cumulativeMass})
making 10 multiplies per cell for applying PPM in one direction for
a Cournat number greater than one. The COSMIC splitting in two dimensions
involves four applications of the one-dimensional PPM. This makes
40 multiplies per cell in total for applying PPM with COSMIC splitting
in two dimensions. In three dimensions, COSMIC splitting requires
12 applications of PPM \citep{LLM96} leading to 120 multiplies.

\subsubsection{Multi-dimensional Scheme}

The number of multiplies involved in the multi-dimensional scheme
includes the number of multiplies to update the higher order advection
and the number of multiplications for each iteration of the linear
equation solver. We will also consider the cost of an explicit version
of the multi-dimensional scheme using an RK2 time-stepping scheme
\citep[eg as used by][]{SW16} which is stable and accurate up to
a Courant number of one for this spatial discretistion and gives very
similar results to the implicit scheme (not shown). 

The explicit version of the dimensionally split scheme uses RK2 or
Heun time-stepping, in which $\phi_{\text{up}}^{\prime}$ on the RHS
of eqn (\ref{eq:CN_1}) is replaced by $\phi_{\text{up}}^{n}$ and
$\phi_{\text{up}}^{n+1}$ on the RHS of eqn (\ref{eq:CN_2}) is replaced
by $\phi_{\text{up}}^{\prime}$. On a logically rectangular two-dimensional
mesh there are 12 cells in each stencil for each face (fig \ref{fig:stencil}).
Each cell has four faces and the interpolation onto each face is used
to calculate the flux between two cells. This leads to 24 multiplies
per cell for each RK2 stage and hence 48 multiplies per cell per time-step.
When using implicit time-stepping, there will be 24 multiplies per
cell for every evaluation of the right hand side of the matrix for
the higher-order correction on first-order upwind. 

We have not explored the sensitivity of the accuracy and stability
to the stencil size and shape in three dimensions. The three dimensional
equivalent of the stencil of quadrilaterals in figure \ref{fig:stencil}
is likely to contain 36 (rather than 12) cells although it may be
possible to omit some corner cells and use a stencil of 20 cells.
For a stencil of 36 cells, the number of multiplies per cell per time-step
would be $36\times3\times2=216$. 

The implicit version of the multi-dimensional scheme (section \ref{sub:implicitMethod})
requires the solution of an asymmetric, diagonally dominant matrix
with three non-zero elements per row for a logically rectangular two-dimensional
mesh. The pre-conditioner is implemented in file \url{DILUPreconditioner.C}
in OpenFOAM 3.0.1 and the solver in file \url{PBiCG.C}. From these
files, we estimate that, for a mesh of quadrilaterals, the solver
will use 24 multiplies (or divides) per cell, per solver iteration,
including pre-conditioning. The average number of iterations of the
preconditioned bi-CG solver per time-step for each of the simulations
is shown in table \ref{tab:iterations} (including the number of solver
iterations in each outer iteration). These simulations use two outer
iterations per time-step when the maximum Courant number is $\le1.1$
(as in equations (\ref{eq:CN_1}) and (\ref{eq:CN_2})) but, for stability,
use four outer iterations per time-step for the larger Courant numbers.
This partly explains the greater number of iterations for larger Courant
numbers. A solver tolerance of $10^{-8}$ is used for each of the
outer solves. Table \ref{tab:iterations} shows that the total number
of iterations per time-step reduces slightly as resolution increases.
The total number of iterations for a complete simulation is reduced
by using larger Courant numbers because the number of iterations per
time-step increases less than linearly with increasing Courant numbers.
In fact, simulations with larger Courant numbers are considerably
cheaper because there are fewer evaluations of the right hand side
of the matrix equation.

\noindent 
\begin{table}
\noindent \begin{centering}
\begin{tabular}{|c||c|c|c|c|c|}
\hline 
Resolution $\rightarrow$ & $120\times60$ & $240\times120$ & $480\times240$ & $960\times480$ & $1920\times960$\tabularnewline
\hline 
Max $c$ $\downarrow$ & \multicolumn{5}{c|}{Iterations per time-step orthogonal/non-orthogonal $\downarrow$}\tabularnewline
\hline 
\hline 
10 & 33.4/39.5 & 31.3/38.0 & 26.0/32.8 & 19.3/25.0 & 13.5/18.6\tabularnewline
\hline 
5 & 22.4/25.7 &  &  &  & \tabularnewline
\hline 
2 & 13.1/14.6 &  &  &  & \tabularnewline
\hline 
1 & 6.0/6.5 & 5.4/5.8 & 4.9/5.3 & 4.6/4.9 & 3.9/4.2\tabularnewline
\hline 
0.5 & 4.9/5.0 &  &  &  & \tabularnewline
\hline 
0.25 & 3.9/3.9 &  &  &  & \tabularnewline
\hline 
\end{tabular}
\par\end{centering}

\caption{Number of solver iterations per time-step of the multi-dimensional
scheme for all the simulations on orthogonal/non-orthogonal meshes.
\label{tab:iterations}}
\end{table}

Combining the number of solver iterations, the number of multiply
operations per solve and the number of multiply operations in calculating
the explicit higher order part of the advection scheme, the total
number of multiply operations per cell per time-step for the multi-dimensional
scheme is shown in table \ref{tab:multiplies} for the non-orthogonal
mesh of $120\times60$ cells. Table \ref{tab:multiplies} also shows
the number of multiplies for using the explicit, RK2 version of the
multi-dimensional scheme and the dimensionally split scheme. 

\begin{table}
\noindent \begin{centering}
\begin{tabular}{|c|c|c|c|}
\hline 
Max $c$ & Implicit & Explicit & Dimensionally\tabularnewline
\hline 
 & \multicolumn{2}{c|}{Multi-dimensional} & Split\tabularnewline
\hline 
\hline 
10 & $39.5\times24+48\times2=1,044$ & - & 40\tabularnewline
\hline 
5 & $25.7\times24+48\times2=712.8$ & - & 40\tabularnewline
\hline 
2 & $14.6\times24+48\times2=446.4$ & - & 40\tabularnewline
\hline 
1 & $6.5\times24+48=204$ & 48 & 40\tabularnewline
\hline 
0.5 & $5.0\times24+48=168$ & 48 & 40\tabularnewline
\hline 
0.25 & $3.9\times24+48=$141.6 & 48 & 40\tabularnewline
\hline 
\end{tabular}
\par\end{centering}

\caption{Total number of multiply operations per cell per time-step for different
Courant numbers for the multi-dimensional scheme, an explicit version
of the multi-dimensional scheme (using RK2 or Heun time-stepping)
and the dimensionally split scheme on two-dimensional, logically rectangular
meshes. The multi-dimensional scheme using RK2 time-stepping is stable
for Courant numbers $\le1$ and gives very similar solutions to the
implicit version. \label{tab:multiplies}}
\end{table}

Table \ref{tab:multiplies} shows that the implicit scheme always
uses more multiply operations but particularly uses more multiplies
for large Courant numbers. The explicit version of the multi-dimensional
scheme always uses 48 multiply operations but is not stable for all
time-steps whereas the dimensionally split scheme (using flux-form
semi-Lagrangian time-stepping) is stable for all Courant numbers (at
this spatial resolution) and always uses the fewest number of multiply
operations per cell per time-step.

There is considerable flexibility in the solver configuration: the
number of outer iterations per time-step determines how frequently
the high-order correction on the right hand side of the matrix equation
is updated, and the solver tolerance per outer iteration could be
modified by using a weaker tolerance on all but the final matrix solve
per time-step. These options have not been explored. It may also be
beneficial to create more non-zero matrix entries rather than having
the higher-order correction entirely a deferred correction on first-order
upwind, but such a change would need to ensure that the matrix remains
diagonally dominant.

\section{Summary and Conclusions\label{sec:concs}}

We examine the errors associated with using a dimensionally split
advection scheme and a multi-dimensional advection scheme on distorted
meshes. The dimensionally split scheme is very accurate on orthogonal
meshes and only loses a little accuracy on highly distorted meshes,
despite a first-order departure point calculation. The multi-dimensional
scheme with implicit time-stepping is less accurate on orthogonal
meshes than the dimensionally split scheme but the accuracy is not
sensitive to mesh distortions and the stability is less sensitive
to Courant number.

The dimensionally split scheme is the piecewise polynomial method
\citep[PPM,][]{CW84} with COSMIC splitting \citep{LLM96} that extends
it to two spatial dimensions. PPM is converges with third-order in
one dimension and COSMIC splitting enables second-order convergence
in two orthogonal directions. PPM is a flux-form semi-Lagrangian scheme
and so can handle large Courant numbers accurately without significant
additional computational cost, with a time-step restriction based
on the deformational Courant number. The second-order accurate multi-dimensional
scheme is split in space and time (method of lines) and uses a cubic
polynomial fit over a stencil of cells for spatial discretisation
and trapezoidal implicit in time to retain stability for large Courant
numbers. We use versions of both schemes without any monotonicity
constraints in order to compare the handling of multi-dimensionality
of the two schemes and order of convergence rather than comparing
the limiters of the two schemes.

Three two-dimensional advection test cases on Cartesian planes are
proposed without the complexities of a spherical domain or multi-panel
meshes but with distorted meshes to mimic the distortions of a cubed-sphere
or terrain following co-ordinates. We therefore propose that these
test cases could be used in the initial testing of advection schemes,
before the generation of meshes on the sphere. The first test case
is an extension of the \citet{LLM96} solid body rotation using a
distorted mesh. The second test case is the established horizontal
advection over orography \citep{SLF+02} using a basic terrain following
mesh in order to maximise the distortions and a version using higher
orography. The third test case is the deformational flow test case
of \citet{LSPT12} adapted to a planar, Cartesian domain and a distorted
mesh. We use the version of this test case with smooth initial conditions
(the sum of two Gaussians) in order to examine order of convergence.

The dimensionally split scheme is extremely accurate on orthogonal
meshes and retains accuracy when long time-steps are used. However
on distorted meshes, particularly at changes in direction such as
those that appear at cube sphere edges or over orography, the split
scheme loses some accuracy. In contrast, the multi-dimensional scheme
is almost entirely insensitive to mesh distortion and asymptotes to
second-order convergence at high resolution. As is expected for implicit
time-stepping, phase errors occur when using long time-steps but the
spatially well resolved features are advected at the correct speed
and the multi-dimensional scheme is always stable.

The matrix solver associated with the implicit time-stepping of the
multi-dimensional scheme means that it is always considerably more
expensive than the split scheme with cost increasing with Courant
number. We haven't investigated sensitivity to solver and pre-conditioner
choices and it may be possible to do  better. But the flux-form semi-Lagrangian
method enables long time-steps without any matrix solutions and this
will always be difficult to beat. It is possible to use the multi-dimensional
scheme with explicit time-stepping such as Runge-Kutta. In that case,
the Courant number is restricted to be less than one but the cost
is similar to the dimensionally split scheme.

The conclusions of this paper are consistent with those of \citet{KNK15}
who found a dimensionally split scheme to be as accurate as a multi-dimensional
scheme on a cubed-sphere mesh with special interpolations so that
the dimensionally split scheme could cope with the cube edges. In
addition, we find that special treatment is not needed at cubed-sphere
edges to maintain accuracy when using dimension splitting.  When dimensionally
split schemes are used on cubed-sphere meshes, they do usually have
special interpolations at cube edges \citep[eg][]{LR96}. However
orography is everywhere and so special treatment over steep orography
would not be practical.

The dimensionally split scheme is hard to beat, providing close to
third-order convergence even in the presence of mesh distortions and
can be very cheaply extended for large Courant numbers.

\bibliographystyle{abbrvnat}
\bibliography{numerics}

\end{document}